\documentclass[12pt]{article}

\usepackage{latexsym}
\usepackage{amsfonts}
\usepackage{amssymb}
\usepackage{amsmath}
\usepackage{mathrsfs}
\input xypic
\usepackage{hyperref}

\newcommand{\be}{\begin{equation}}
\newcommand{\ee}{\end{equation}}
\newcommand{\bea}{\begin{eqnarray}}
\newcommand{\eea}{\end{eqnarray}}
\newcommand{\bean}{\begin{eqnarray*}}
\newcommand{\eean}{\end{eqnarray*}}
\newcommand{\brray}{\begin{array}}
\newcommand{\erray}{\end{array}}

\newtheorem{dfn}{Definition}[section]
\newtheorem{thm}[dfn]{Theorem}
\newtheorem{lmma}[dfn]{Lemma}
\newtheorem{ppsn}[dfn]{Proposition}
\newtheorem{crlre}[dfn]{Corollary}
\newtheorem{xmpl}[dfn]{Example}
\newtheorem{rmrk}[dfn]{Remark}

\newcommand{\bdfn}{\begin{dfn}\rm}
\newcommand{\bthm}{\begin{thm}}
\newcommand{\blmma}{\begin{lmma}}
\newcommand{\bppsn}{\begin{ppsn}}
\newcommand{\bcrlre}{\begin{crlre}}
\newcommand{\bxmpl}{\begin{xmpl}}
\newcommand{\brmrk}{\begin{rmrk}\rm}

\newcommand{\edfn}{\end{dfn}}
\newcommand{\ethm}{\end{thm}}
\newcommand{\elmma}{\end{lmma}}
\newcommand{\eppsn}{\end{ppsn}}
\newcommand{\ecrlre}{\end{crlre}}
\newcommand{\exmpl}{\end{xmpl}}
\newcommand{\ermrk}{\end{rmrk}}

\newcommand{\cla}{\mathcal{A}}

\newcommand{\clh}{\mathcal{H}}

\newcommand{\clg}{\mathcal{G}}

\author{S. Sundar }
\title{ $C^{*}$-algebras associated to  topological Ore semigroups }

\begin{document}
\maketitle 
\begin{abstract}
Let $G$ be a locally compact group and $P \subset G$ be a closed Ore semigroup containing the identity element. Let $V: P \to B(\clh)$ be an anti-homomorphism such that for every $a \in P$, $V_{a}$ is an isometry and the final projections of $\{V_{a}: a \in P\}$ commute. We study the $C^{*}$-algebra generated by $\{\int f(a)V_{a} da: f \in L^{1}(P)\}$. We show that there exists a groupoid $C^{*}$-algebra which is universal for isometric representations with commuting range projections.  \end{abstract}

\noindent {\bf AMS Classification No. :} {Primary 22A22; Secondary 54H20, 43A65, 46L55.}  

\section{Introduction}

It is fair to say that $C^{*}$-algebras of groups and their crossed products are the most studied $C^{*}$-algebras in the theory of operator algebras. Several authors have tried to study $C^{*}$-algebras associated to  semigroups. For example, the Toeplitz algebra is the $C^{*}$-algebra associated to the additive semigroup $\mathbb{N}$. Recently, the theory of semigroup $C^{*}$-algebras have received renewed attention. See for example \cite{Cuntz}, \cite{Li-semigroup}, \cite{Li13} and the references therein.  The notion of  crossed product by semigroups has also been studied by several authors most notably by Murphy in \cite{Murphy91}, \cite{Murphy94}, \cite{Murphy96} and by Exel  in \cite{Exel_endo}. However much of the literature focusses on discrete semigroups. In the topological direction, upto the author's knowledege, the only example studied is the Wiener-Hopf $C^{*}$-algebra. This was studied from the groupoid point of view first in \cite{Renault_Muhly} and then successively by Nica in \cite{Nica_WienerHopf},  \cite{Nica90}
 and  Hilgert and Neeb in \cite{Hilgert_Neeb}.

 Let $G$ be a second countable locally compact group and $P \subset G$ be a closed semigroup containing the identity element. We assume that $Int(P)$ is dense in $P$ and $PP^{-1}=G$.  Let $V: P \to B(\clh)$ be an isometric representation on a Hilbert space $\clh$ i.e. for $a \in P$, $V_{a}$ is an isometry and $V_{a}V_{b}=V_{ba}$. For $f \in L^{1}(P)$, let 
 \[
 W_{f}:= \int_{a \in P} f(a)V_{a} da.
  \]
The semigroup $C^{*}$-algebra or the Wiener-Hopf algebra, denoted $\mathcal{W}_{V}(P,G)$, associated to the representation $V$ is the $C^{*}$-algebra generated by $\{W_{f}:f \in L^{1}(P)\}$.  If we consider the compression of the right regular representation of $G$ on $L^{2}(G)$ onto $L^{2}(P)$, then one obtains the usual Wiener-Hopf algebra studied in \cite{Renault_Muhly}. In general, it is much difficult to understand the structure of $\mathcal{W}_{V}(P,G)$. However if we assume that the final projections $\{E_{a}:=V_{a}V_{a}^{*}: a \in P\}$ form a commuting family of projections then one can do better. Without this commutative assumption, the situation becomes much complicated even for the simplest case of $P:=\mathbb{N} \times \mathbb{N}$ as is illustrated by Murphy in \cite{Murphy}. The results obtained and the organisation of the paper are described below.

From now on, we assume that the range projections commute. For $g =ab^{-1} \in G$, let $W_{g}:= V_{b}^{*}V_{a}$ and $ E_{g}$ be the final space of $W_{g}$. It is shown in Section 3, that $W_{g}$ is well-defined and $\{E_{g}:g \in G\}$ forms a commuting family of projections. For $f \in L^{1}(G)$, let $W_{f}:=\int f(g)W_{g} dg$. It is not difficult to show that $\mathcal{W}_{V}(P,G)$ is generated by $\{W_{f}: f \in L^{1}(G)\}$. Let $\Omega$ be the spectrum of the commutative $C^{*}$-algebra generated by $\{\int f(g)E_{g} dg: f \in L^{1}(G)\}$.   The map $C(\Omega) \rtimes P \ni (T,a) \to V_{a}^{*}TV_{a} \in C(\Omega)$ provides an action of $P$ on $\Omega$. In Section 4 and 5, we show that this action is injective. Let \[\mathcal{G}:= \Omega \rtimes P:=\{(x,ab^{-1},y) \in \Omega \times G \times \Omega: xa=yb\}\] be the Deaconu-Renault groupoid where the groupoid operations are given by 
\begin{align*}
 (x,g,y)(y,h,z)&=(x,gh,z) \\
 (x,g,y)^{-1}&=(y,g^{-1},x).
\end{align*}
 For $f \in C_{c}(G)$, let $\widetilde{f} \in C_{c}(\mathcal{G})$ be defined by $\widetilde{f}(x,g,y)=f(g)$.   
 We apply the results of \cite{Jean_Sundar} to show that $\Omega \rtimes P$ has a Haar system.  We also show that there exists a surjective representation $\lambda: C^{*}(\mathcal{G}) \to \mathcal{W}_{V}(P,G)$ such that for $f \in C_{c}(G)$,
 \[
 \lambda(\widetilde{f})=\int f(g)\Delta(g)^{-\frac{1}{2}} W_{g^{-1}} dg.
 \]
Here $\Delta$ denotes the modular function of the group. This is achieved in Sections 4-6. For the Wiener-Hopf representation, the groupoid $\Omega \rtimes P$ is the groupoid considered in \cite{Renault_Muhly}. 

We  show in Section 7, that there exists a universal space $\Omega_{u}$ on which $P$ acts such that if $V:P \to B(\clh)$ is an isometric representation with commuting range projections then there exists a representation $\lambda:C^{*}(\Omega_{u} \rtimes P) \to B(\clh)$ such that for $f \in C_{c}(G)$, 
\[
\lambda(\widetilde{f})=\int f(g)\Delta(g)^{-\frac{1}{2}}W_{g^{-1}} dg.
\]

\section{Preliminaries}
For the convenience of the reader, we recall the essential facts from \cite{Jean_Sundar} that we need in this paper. The proofs can be found in \cite{Jean_Sundar}. Throughout this paper,  $G$  stands for a second countable, locally compact topological group and $P \subset G$ for a closed subsemigroup containing the identity element $e$. We also assume the following. 
\begin{enumerate}
 \item[(C1)] The group $G=PP^{-1}$, and
  \item[(C2)] the interior of $P$ in $G$, denoted $Int(P)$, is dense in $P$.
\end{enumerate}
Semigroups for which $(C1)$ is satisfied are called Ore semigroups. In this paper, we consider only semigroups with identity for which $(C1)$ and $(C2)$ are satisfied.

Let $X$ be a compact Hausdorff space. A right action of $P$ on $X$ is a continuous map $ X \times P \ni (x,a) \to xa \in X$ such that $xe=x$ and $(xa)b=x(ab)$ for $x \in X$ and $a,b \in P$. Moreover we assume that the action is injective i.e. for every $a \in P$, the map $X \ni x \to xa \in X$ is injective. Let $X$ be a compact Hausdorff space on which $P$ acts on the right injectively. Then the semi-direct product groupoid $X \rtimes P$ is defined as follows: \[ X \rtimes P: = \{(x,g,y) \in X \times G \times X: \exists ~a,b \in P,\textrm{~such that~}g=ab^{-1}, xa=yb \}. \]
The groupoid multiplication and the inversion are given by 
\begin{align*}
 (x,g,y)(y,h,z)&=(x,gh,z), \\
 (x,g,y)^{-1}&=(y,g^{-1},x).
\end{align*}

 The map $ X \rtimes P \ni (x,g,y) \to (x,g) \in X \times G$ is injective. Thus  $X \rtimes P$ can be considered a subset of $X \times G$ which we do from now. Moreover $X \rtimes P$ is a closed subset of $X \times G$ and when $X \rtimes P$ is given the subspace topology, the groupoid $X \rtimes P$ becomes a topological groupoid. We denote the range and source maps by $r$ and $s$ respectively.

For $x \in X$, let $Q_{x}:= \{g \in G: (x,g) \in \mathcal{G}\}$. Then $r^{-1}(x)= \{x\} \times Q_{x}$. Note that for $x \in X$, $Q_{x}.P \subset Q_{x}$ and $Q_{x}$ is closed. By Lemma 4.1 of \cite{Jean_Sundar},  for $x \in X$, $Int(Q_{x})$ is dense in $Q_{x}$ and the boundary of $Q_{x}$ has measure zero.

For $x \in X$, let $\lambda^{x}$ be the measure on $\mathcal{G}$ defined as follows: For $f \in C_{c}(\mathcal{G})$,
\[
 \int f d\lambda^{x} = \int f(x,g)1_{Q_{x}}(g)dg.
\]
Here $dg$ denotes the left Haar measure on $G$. In \cite{Jean_Sundar}, it is shown that the groupoid $\mathcal{G}:=X \rtimes P$ admits a Haar system if and only if the map $X \times Int(P) \ni (x,a) \to xa \in X$ is open. In this case,  the measures $(\lambda^{x})_{x \in X}$ form  a Haar system.  We will  use only this Haar system if $X \rtimes P$ admits one.

Suppose that $\mathcal{G}:=X \rtimes P$ admits a Haar system. Then the action of $P$ on $X$ can be dilated to an action of $G$. That is there exists a locally compact  Hausdorff space $Y$ on which $G$ acts on the right and a continuous $P$-equivariant injection $i:X \to Y$ such that 
\begin{enumerate}
 \item the set $X_{0}:=i(X)Int(P)$ is open in $Y$, and
 \item $Y=\bigcap_{a\in P}i(X)a^{-1}=\bigcap_{a \in Int(P)}X_{0}a^{-1}$.
\end{enumerate}
Moreover the space $Y$ is unique up to a $G$-equivariant homeomorphism. We will identiy $X$ as a subspace of $Y$ via the injection $i$ and will suppress the notation $i$. Also the groupoid $\mathcal{G}$ is isomorphic to the reduction $(Y \rtimes G)|_{X}$. 
With this notation, note that for $x \in X$, $1_{Q_{x}}(g)=1_{X}(xg)$. Also we leave it to the reader to check that $1_{Int(Q_{x})}(g)=1_{X_{0}}(xg)$.

For $f \in C_{c}(G)$, let $\widetilde{f} \in C_{c}(\mathcal{G})$ be defined by $\widetilde{f}(x,g)=f(g)$. 
We  also need the following proposition. The proof is a line by line imitation of that of Proposition 3.5 of \cite{Renault_Muhly}. Hence we omit the proof. See also \cite{Jean_Sundar} for some remarks concerning the proof. 

Let $\mathcal{G}:=X \rtimes P$ and assume that it has a Haar system. 

\begin{ppsn}
\label{density}
For $f \in C_{c}(G)$, let $\hat{f} \in C(X)$ be defined by $\hat{f}(x)=\int f(g)1_{X}(xg)dg$. Suppose that family $\{\hat{f}: f \in C_{c}(G)\}$ separates points of $X$. Then  the $*$-algebra generated by $\{\widetilde{f}: f \in C_{c}(G)\}$  is dense in $C_{c}(\mathcal{G})$ where $C_{c}(\mathcal{G})$ is given the inductive limit topology. As a consequence, $C^{*}(\mathcal{G})$ is generated by $\{\widetilde{f}:f \in C_{c}(G)\}$.
 \end{ppsn}

\section{Isometric representations with commuting range projections}

\begin{dfn}
 A map $V:P \to B(\clh)$ is called an  isometric representation of $P$ on the Hilbert space $\clh$ if 
\begin{enumerate}
 \item [(1)] the maps $P \ni a \to V_{a}$ and $P \ni a \to V_{a}^{*}$ are strongly continuous,
 \item [(2)] for $a \in P$, $V_{a}$ is an isometry, and
 \item [(3)] for $a,b \in P$, $V_{a}V_{b}=V_{ba}$.
\end{enumerate}
 For $a \in P$, let $E_{a}:=V_{a}V_{a}^{*}$.  If $\{E_{a}: a\in P\}$ is a commuting family of projections, we say that $V$ has commuting range projections.
\end{dfn}

In the next example, we recall the Wiener-Hopf representation or the regular representation. The $C^{*}$-algebra associated to the Wiener-Hopf representation has been studied by several authors. See the papers \cite{Renault_Muhly}, \cite{Hilgert_Neeb} and the references therein.

\begin{xmpl}
Consider the Hilbert space $L^{2}(G)$ and consider $L^{2}(P)$ as a closed subspace of $L^{2}(G)$. 
For $\xi \in L^{2}(P)$ and $a \in P$, let $V_{a}(\xi)$ be defined as follows: 
\begin{displaymath}
\begin{array}{lll}
V_{a}(\xi)(x):=& \left\{\begin{array}{lll}
                               \xi(xa^{-1})\Delta(a)^{\frac{-1}{2}} &~if~x \in Pa \\
                                 0 & ~if~x \notin Pa.
                                 \end{array} \right. 
\end{array}
\end{displaymath}
Here $\Delta$ denotes the modular function of the group $G$. Then the map $a \in P \to V_{a} \in B(L^{2}(P))$ is an isometric representation with commuting range projections. Note that for $a \in P$, the range of $V_{a}$ is $L^{2}(Pa)$.
\end{xmpl}

Till the end of Section 7, we fix an isometric representation $V:P\to B(\clh)$ with commuting range projections. 
For $g=ab^{-1}$, let $W_{g}:=V_{b}^{*}V_{a}$ and let $E_{g}:=W_{g}W_{g}^{*}$. First we show that that $W_{g}$ is well defined and is a partial isometry.

\begin{ppsn}
\label{commuting projections}
 Let $V:P \to B(\clh)$ be an isometric representation with commuting range projections. 
\begin{enumerate}
 \item[(1)] For $g \in G$, $W_{g}$ is well defined and is a partial isometry.
 \item[(2)] The family $\{E_{g}:g \in G \}$ forms a commuting family of projections.
  \item[(3)] If $g_{1}g_{2}^{-1} \in P$, then $E_{g_{1}} \leq E_{g_{2}}$.
  \item[(4)] The map $G \in g \to W_{g} \in B(\clh)$ is strongly continuous.
  \item[(5)] For $g, h \in G$, $W_{g}W_{h}=E_{g}W_{hg}$.
\end{enumerate}
\end{ppsn}
\textit{Proof.} Suppose   $g=a_{1}b_{1}^{-1}=a_{2}b_{2}^{-1}$. Then $a_{1}^{-1}a_{2}=b_{1}^{-1}b_{2}$. Since $PP^{-1}=G$, there exists  $\alpha_{1},\alpha_{2} \in P$ such that $a_{1}^{-1}a_{2}=\alpha_{1} \alpha_{2}^{-1}=b_{1}^{-1}b_{2}$. Then $a_{1}\alpha_{1}=a_{2}\alpha_{2}$ and $b_{1}\alpha_{1}=b_{2}\alpha_{2}$.

Now observe that
\begin{align*}
 V_{b_{1}}^{*}V_{a_{1}}  &= V_{b_{1}}^{*}V_{\alpha_{1}}^{*}V_{\alpha_{1}}V_{a_{1}} \\
                         & =V_{b_{1}\alpha_{1}}^{*}V_{a_{1}\alpha_{1}} \\
                         &= V_{b_{2}\alpha_{2}}^{*}V_{a_{2}\alpha_{2}} \\
                         & =V_{b_{2}}^{*}V_{\alpha_{2}}^{*}V_{\alpha_{2}}V_{a_{2}} \\
                          & =V_{b_{2}}^{*}V_{a_{2}}.
\end{align*}
This proves that $W_{g}$ is well defined. Let $E_{g}:=W_{g}W_{g}^{*}$. If $g=ab^{-1}$, then $E_{g}=V_{b}^{*}E_{a}V_{b}$ which is self adjoint. Now note that 
\begin{align*}
 E_{g}^{2} &=V_{b}^{*}E_{a}V_{b}V_{b}^{*}E_{a}V_{b}\\
           & = V_{b}^{*}E_{a}E_{b}E_{a}V_{b} \\
           &=V_{b}^{*}E_{b}E_{a}^{2}V_{b} ~~~\textrm{( Since $E_{a}$ and $E_{b}$ commute)} \\
           &=V_{b}^{*}E_{a}V_{b} \\
           &=E_{g}.
\end{align*}
Thus $E_{g}$ is a projection. This proves $(1)$.

Let $g_{1},g_{2} \in G$ be given. Write $g_{1}=a_{1}b_{1}^{-1}$ and $g_{2}=a_{2}b_{2}^{-1}$ with $a_{i},b_{i} \in P$. Choose $\alpha_{1}, \alpha_{2} \in P$ such that $b_{1}\alpha_{1}=b_{2}\alpha_{2}$. Let $a_{i}^{'}=a_{i}\alpha_{i}$ and $b_{i}^{'}=b_{i}\alpha_{i}$ for $i=1,2$. Then $g_{i}=a_{i}^{'}(b_{i}^{'})^{-1}$ for $i=1,2$. But now $b_{1}^{'}=b_{2}^{'}$.
Thus 
\begin{align*}
 E_{g_{1}}E_{g_{2}}&=V_{b_{1}^{'}}^{*}E_{a_{1}^{'}}V_{b_{1}^{'}}V_{b_{1}^{'}}^{*}E_{a_{2}^{'}}V_{b_{1}^{'}} \\
                   & = V_{b_{1}^{'}}^{*}E_{a_{2}^{'}}V_{b_{1}^{'}}V_{b_{1}^{'}}^{*}E_{a_{1}^{'}}V_{b_{1}^{'}} \\
                   & = E_{g_{2}}E_{g_{1}}.
\end{align*}
This proves $(2)$.

Suppose $g_{1}g_{2}^{-1}=a$ for some $a \in P$. Write $g_{2}=bc^{-1}$. Then $g_{1}=(ab)c^{-1}$. Then 
\begin{align*}
E_{g_{1}}&=V_{c}^{*}V_{ab}V_{ab}^{*}V_{c} \\
         &= V_{c}^{*}V_{b}(V_{a}V_{a}^{*})V_{b}^{*}V_{c} \\
         & \leq V_{c}^{*}V_{b}V_{b}^{*}V_{c} \\
         & \leq E_{g_{2}}.
\end{align*}
This proves $(3)$.

Note that the map $Int(P) \times Int(P) \ni (a,b) \to ab^{-1} \in G$ is surjective and open. Thus $G$ is the quotient of $Int(P) \times Int(P)$. Since   multiplication is strongly continuous on the unit ball of $B(\clh)$, it follows that the map $Int(P) \times Int(P) \ni (a,b) \to V_{b}^{*}V_{a} \in B(\clh)$ is strongly continuous. As a consequence, it follows that the map $G \ni g \to W_{g} \in B(\clh)$ is strongly continuous. This proves $(4)$. 

Let $g,h \in G$ be given. Write $g=ab^{-1}$ and $h=cd^{-1}$ with $a,b \in P$. Choose $\alpha,\beta \in P$ such that $d\beta=a\alpha$. Now note that $g=(a\alpha)(b\alpha)^{-1}$ and $h=(c\beta)(d\beta)^{-1}$. Thus we can write $g$ and $h$ as $g=a_{1}b_{1}^{-1}$ and $h=c_1a_{1}^{-1}$. Now calculate as follows
\begin{align*}
W_{g}W_{h}&=V_{b_1}^{*}V_{a_1}V_{a_1}^{*}V_{c_1} \\
          &=V_{b_1}^{*}E_{a_1}V_{b_{1}}V_{b_{1}}^{*}V_{c_1} (\textrm{~Since $E_{a_1}$ and $E_{b_1}$ commute}) \\
          &=E_{g}W_{hg}
\end{align*}

This completes the proof. \hfill $\Box$

For $f \in L^{1}(G)$, the ``Wiener-Hopf'' operator with symbol $f$ is defined as \[ W_{f}:=\int f(g)W_{g}dg. \]
We want to describe the $C^{*}$-algebra, denoted $\mathcal{W}_{V}(P,G)$, generated by  $\{W_{f}: f \in L^{1}(G)\}$. We suppress the subscript $V$ and simply denote $\mathcal{W}_{V}(P,G)$ by $\mathcal{W}(P,G)$ ( atleast till the end of Section 7.)

\begin{rmrk}
One can show that $\mathcal{W}(P,G)$ is generated by $\{\int f(a)V_{a}da: f \in L^{1}(P)\}$. The proof is similar to that of Proposition 2.2 of \cite{Jean_Sundar}. Hence we omit the proof.
\end{rmrk}

First we consider  a related commutative $C^{*}$-algebra. Note that by definition, for $g \in G$, $W_{g^{-1}}=W_{g}^{*}$. Moreover the map $G \ni g \to E_{g}=W_{g}W_{g}^{*}$ is strongly continuous. For $f \in L^{1}(G)$,  let \[
                           E_{f}:= \int f(g)E_{g}dg.\]
                         
Let \[
     \mathcal{A}:= C^{*}\{ E_{f}: f \in L^{1}(G) \}.
    \]
Since $\{E_{g}:g \in G\}$ forms a commuting family of projections, it follows that $\cla$ is a commutative $C^{*}$-subalgebra of $B(\clh)$. Note that $E_{g}=1$ if $g \in P^{-1}$. If $f \in L^{1}(P^{-1})$, then $E_{f} = \int f(g)dg$.   Thus, it follows that $\mathcal{A}$ is a commutative unital $C^{*}$ subalgebra of $B(\clh)$. Denote the spectrum of $\mathcal{A}$ by $\Omega$.

Let $G_{n}:= \underbrace{G \times G \times \cdots \times G}_{\textrm{$n$ times}}$. For $f \in C_{c}(G_{n})$, let 
\[
 E_{f}:= \int f(g_{1},g_{2},\cdots,g_{n})E_{g_{1}}E_{g_{2}}\cdots E_{g_{n}}dg_{1}dg_{2}\cdots dg_{n}.
\]
Let $\displaystyle \widetilde{\mathcal{A}}:= \bigcup_{n=1}^{\infty}\{E_{f}: f \in C_{c}(G_{n})\}$. Then $\widetilde{\mathcal{A}}$ forms a dense unital $*$-subalgebra of $\cla$. Also note that for every $n$, the map $C_{c}(G_{n}) \ni f \to E_{f} \in \mathcal{A}$ is continuous when $C_{c}(G_{n})$ is given the inductive limit topology and $\mathcal{A}$ is given the norm topology.

For $T \in B(\clh)$ and $a \in P$, let $\alpha_{a}(T)=V_{a}^{*}TV_{a}$. Clearly $\alpha_{e}=id$ and $\alpha_{a}\alpha_{b}=\alpha_{ab}$.  

 Observe that  $\alpha_{a}(V_{b}^{*}E_{c}V_{b})=V_{a}^{*}V_{b}^{*}E_{c}V_{b}V_{a}=V_{ab}^{*}E_{c}V_{ab}$. Thus  
$\alpha_{a}(E_{g})=E_{ga^{-1}}$ for $g \in G$.  Since the final projection $V_{a}V_{a}^{*}$ commutes with $E_{g}$ for every $g \in G$, it follows that $\alpha_{a}(E_{g_{1}}E_{g_{2}})=\alpha_{a}(E_{g_{1}})\alpha_{a}(E_{g_{2}})$.

\begin{ppsn}
\label{The action}
 For $a \in P$, $\alpha_{a}$ leaves $\mathcal{A}$ invariant and the map $\alpha_{a}:\mathcal{A} \to \mathcal{A}$ is a unital $*$-homomorphism. Moreover for $T \in \mathcal{A}$, the map $P \ni a \to \alpha_{a}(T) \in \mathcal{A}$ is norm continuous.
\end{ppsn}
\textit{Proof.}
For $a \in P$ and $f \in C_{c}(G_{n})$, let $\widetilde{f}_{a} \in C_{c}(G_{n})$ be defined by \[
                                                                                                \widetilde{f}_{a}(g_{1},g_{2},\cdots,g_{n})=\Delta(a)^{n}f(g_{1}a,g_{2}a,\cdots,g_{n}a).
                                                                                               \]
Then for $f \in C_{c}(G_{n})$, the map $P \ni a \to \widetilde{f}_{a} \in C_{c}(G_{n})$ is continuous if $C_{c}(G_{n})$ is given the inductive limit topology.

Let $a \in P$ and $f \in C_{c}(G_{n})$ be given. Then \begin{align*}
                                                       \alpha_{a}(E_{f}) &= \int f(g_{1},g_{2},\cdots,g_{n})\alpha_{a}(E_{g_{1}}E_{g_{2}}\cdots E_{g_{n}})dg_{1}dg_{2}\cdots dg_{n} \\
                                                                         & = \int f(g_{1},g_{2},\cdots,g_{n})E_{g_{1}a^{-1}}E_{g_{2}a^{-1}}\dots E_{g_{n}a^{-1}}dg_{1}dg_{2} \cdots dg_{n} \\
                                                                        & = \int f(g_{1}a,g_{2}a,\cdots,g_{n}a)\Delta(a)^{n}E_{g_{1}}E_{g_{2}}\cdots E_{g_{n}}dg_{1}dg_{2}\cdots dg_{n} \\
                                                                        & = E_{\widetilde{f}_{a}}
                                                      \end{align*}
 Thus $\alpha_{a}$ leaves $\widetilde{\mathcal{A}}$ invariant. Since $\widetilde{\mathcal{A}}$ is dense in $\mathcal{A}$ and $\alpha_{a}$ is bounded, it follows that $\alpha_{a}$ leaves $\mathcal{A}$ invariant.

Observe that if $E_{a}=V_{a}V_{a}^{*}$ commutes with $T, S \in B(\clh)$ then $\alpha_{a}(TS)=\alpha_{a}(T)\alpha_{a}(S)$. By Proposition \ref{commuting projections}, if follows that  $E_{a}$ commutes with $E_{f}$ for $f \in C_{c}(G_{n})$. Thus $E_{a}$ commutes with every element of $\cla$. Hence $\alpha_{a}:\mathcal{A} \to \mathcal{A}$ is multiplicative. Clearly $\alpha_{a}$ is unital and $*$-preserving.

For $a \in P$, $\alpha_{a}:\mathcal{A} \to \mathcal{A}$ is contractive. Thus it is enough to show that  for $T \in \widetilde{\mathcal{A}}$, the map $P \ni a \to \alpha_{a}(T) \in \mathcal{A}$ is continuous. Let $T=E_{f}$ for some $f \in C_{c}(G_{n})$. Then $\alpha_{a}(T)=E_{\widetilde{f}_{a}}$. Hence the map $P \ni a \to \alpha_{a}(T)=E_{\widetilde{f}_{a}}$ is continuous as it is the composite of the continuous maps $P \ni a \to \widetilde{f}_{a} \in C_{c}(G_{n})$ and $C_{c}(G_{n}) \ni h \to E_{h}$ where $C_{c}(G_{n})$ is given the inductive limit topology. This completes the proof. \hfill $\Box$

Since $\mathcal{A}=C(\Omega)$, it follows that for every $a \in P$, there exists $\phi_{a}:\Omega \to \Omega$ such that $F\circ \phi_{a}=\alpha_{a}(F)$ for $F \in C(\Omega)$. The condition $\alpha_{a}\alpha_{b}=\alpha_{ab}$ translates to $\phi_{a}\phi_{b}=\phi_{ba}$ for $a,b \in P$. Also $\phi_{e}=id$. Thus the map $\Omega \times P \ni (x,a) \to \phi_{a}(x) \in \Omega$ defines a right action of $P$ on $\Omega$. We henceforth write $\phi_{a}(x)$ as $xa$ for $x \in \Omega$ and $a \in P$.

We claim that the map $\Omega \times P \ni (x,a) \to xa \in \Omega$ is continuous. Suppose $(x_{n}) \to x$ and $(a_{n}) \to a$. Let $ F \in C(\Omega)$. By Proposition \ref{The action}, it follows that $\alpha_{a_{n}}(F)$ converges uniformly to $\alpha_{a}(F)$. Since the convergence is uniform, it follows that $\alpha_{a_{n}}(F)(x_n)$ converges to $\alpha_{a}(F)(x)$. In other words, for every $F \in C(\Omega)$, $F(x_{n}a_{n})$ converges to $F(xa)$. Hence $x_{n}a_{n}$ converges to $xa$.

The goal of this paper is to prove the following statements.

\begin{enumerate}
 \item[(1)] The right action of $P$ on $\Omega$ is injective.
 \item[(2)] The semidirect product groupoid $\mathcal{G}:=\Omega \rtimes P$ has a Haar system.
 \item[(3)] For $f \in C_{c}(G)$, let $\widetilde{f} \in C_{c}(\mathcal{G})$ be defined by $\widetilde{f}(x,g)=f(g)$ for $(x,g) \in \mathcal{G}$. There exists a surjective $*$-homomorphism $\pi:C^{*}(\mathcal{G}) \to \mathcal{W}(P,G)$ such that $\pi(\widetilde{f})=\int \Delta(g)^{-\frac{1}{2}}f(g)W_{g^{-1}}dg$ for $f \in C_{c}(G)$.
\end{enumerate}
To prove the above statements, we need a better description of $\Omega$ which forms the content of the next section. We end this section with a lemma which is useful in showing that $\Omega \rtimes P$ has a Haar system.

\begin{lmma}
 \label{Transfer operator}
Let $f \in C_{c}(G)$ be such that $supp(f) \subset Int(P)$. Then for $T \in \mathcal{A}$, the integral $\displaystyle \int_{a \in P} f(a)V_{a}TV_{a}^{*}da \in \mathcal{A}$. 
\end{lmma}
\textit{Proof.} It is enough to prove the statement for $T \in \widetilde{\mathcal{A}}$. Let $T=E_{\phi}$ for some $\phi \in C_{c}(G_{n})$. For $a \in P$ and $g \in G$, $V_{a}^{*}E_{ga}V_{a}=E_{g}$. Hence $V_{a}E_{g}V_{a}^{*}=E_{ga}E_{a}$. Now calculate as follows to find that
\begin{align*}
 \int_{a \in P} &f(a)V_{a}E_{\phi}V_{a}^{*}da& \\
&= \int_{a \in P} f(a)\phi(g_{1},g_{2},\cdots,g_{n})V_{a}E_{g_{1}}E_{g_{2}}\cdots E_{g_{n}}V_{a}^{*}da ~dg_{1}~ dg_{2} \cdots dg_{n} \\
       & = \int_{a \in P}f(a)\phi(g_{1},g_{2},\cdots,g_{n})E_{a}E_{g_{1}a}E_{g_{2}a}\cdots E_{g_{n}a}~da~ dg_{1}dg_{2} \cdots dg_{n} \\
                 & = \int_{a \in P}\Delta(a)^{-n}f(a)\phi(g_{1}a^{-1},g_{2}a^{-1},\cdots,g_{n}a^{-1})E_{a}E_{g_{1}}E_{g_{2}}\cdots E_{g_{n}}~da~ dg_{1}~dg_{2}\cdots dg_{n} \\
                 & = E_{\psi} \in \widetilde{\mathcal{A}}
\end{align*}
where $\psi \in C_{c}(G_{n+1})$ is given by \[\psi(g,g_{1},g_{2},\cdots,g_{n})=\Delta(g)^{-n}f(g)\phi(g_{1}g^{-1},g_{2}g^{-1},\cdots,g_{n}g^{-1}).\]
This completes the proof. \hfill $\Box$

\section{What is  $\Omega$ ? }

We first discuss the case when $G$ is discrete. The discrete semigroup $C^{*}$-algebras are analysed in great detail  in the papers \cite{Li-semigroup} and \cite{Li13}. Neverthless we discuss this case in the form that we need. This also motivates the topological case.

Let $G$ be a discrete group and $P \subset G$ be a semigroup such that $e \in P$ and $PP^{-1}=G$. In this case, the Wiener-Hopf $C^{*}$-algebra $\mathcal{W}_{V}(P,G)$ is simply the $C^{*}$-algebra generated by $\{W_{g}:g \in G\}$ and the commutative $C^{*}$-algebra $\mathcal{A}$ is the $C^{*}$-algebra generated by $\{E_{g}:g \in G\}$.

Let $\chi$ be a character of $\mathcal{A}$. Let us define the support of $\chi$, denoted  $A_{\chi}$,  as
\[
 A_{\chi}:=\{g \in G: \chi(E_{g})=1\}.
\]
 Condition $(3)$ of Proposition \ref{commuting projections} implies that $P^{-1}A_{\chi} \subset A_{\chi}$. Since $E_{g}=1$ if $g \in P^{-1}$, it follows that $P^{-1} \subset A_{\chi}$. 

Let $\mathcal{P}(G)$ be the power set  of $G$. Identify $\mathcal{P}(G)$ with  $\{0,1\}^{G}$, via the map $\mathcal{P}(G) \ni A \to 1_{A} \in \{0,1\}^{G}$, and endow it with the product topology. The group $G$ acts on $\mathcal{P}(G)$. The right action is  given by : For $g \in G$ and $A \in \mathcal{P}(G)$, $Ag:= \{ag: a \in A\}$. Clearly the map $ \Omega \ni \chi \to A_{\chi} \in \{0,1\}^{G}$ is continuous, injective and hence an embedding. We leave it to the reader to check that the above map is $P$-equivariant. From now, we view $\Omega$ as a subset of $\{0,1\}^{G}$.

\begin{ppsn}
\label{discrete case}
We have the following.
\begin{enumerate}
 \item[(1)] For $A \in \Omega$ and $a \in P$, $Aa^{-1} \in \Omega$ if and only if $a \in A$.
  \item[(2)] For $A \in \Omega$ and $g \in G$, $Ag \in \Omega$ if and only if $g^{-1} \in A$.
 \item[(3)] The action $\Omega \times P \to \Omega$ is open.
\end{enumerate}
\end{ppsn}
\textit{Proof.} Let $A \in \Omega$ and $a \in P$ be given. Suppose $B:=Aa^{-1} \in \Omega$. Since $e \in B$, it follows that $a \in A$. Now suppose $a \in A$. Let $\chi$ be the character corresponding to $A$. Since $V_{a}^{*}E_{ga}V_{a}=E_{g}$, it follows that $V_{a}E_{g}V_{a}^{*}=E_{a}E_{ga}$. Thus the homomorphism $B(\clh) \ni T \to V_{a}TV_{a}^{*} \in B(\clh)$ leaves $\mathcal{A}$ invariant. 

Let $\widetilde{\chi}$ be the character on $\mathcal{A}$  defined by 
$\widetilde{\chi}(T)=\chi(V_{a}TV_{a}^{*})$. Since $\chi(E_a)=1$, it follows that $\widetilde{\chi}$ is non-zero. Observe that for $T \in \mathcal{A}$, 
\begin{align*}
 (\widetilde{\chi}a)(T)&=\widetilde{\chi}(V_{a}^{*}TV_{a}) \\
                       &=\chi(V_{a}V_{a}^{*}TV_{a}V_{a}^{*})\\
                       &= \chi(E_{a})\chi(T)\chi(E_{a}) \\
                       &= 1_{A}(a)\chi(T)1_{A}(a) \\
                       & = \chi(T).
\end{align*}
Thus $\widetilde{\chi}a=\chi$. Let $B$ be the support of $\widetilde{\chi}$. Then $A=Ba$. Thus $Aa^{-1} \in \Omega$. This proves $(1)$.

Now let $A \in \Omega$ and $g=ab^{-1} \in G$. Suppose $g^{-1}=ba^{-1} \in A$. Then $b \in Aa \in \Omega$. By $(1)$,
it follows that $Ag=Aab^{-1} \in \Omega$. Now suppose $Ag \in \Omega$. Then $A=(Ag)g^{-1}$. Since $e \in Ag$, it follows that $g^{-1} \in A$. This proves $(2)$.

When $G$ is discrete, $Int(P)=P$ and $\Omega P=\Omega$. Thus,
by Theorem 4.3 of \cite{Jean_Sundar}, to prove that the action $\Omega \times P \to \Omega$ is open, it is enough to show that $\Omega a$ is open in $\Omega$ for every $a \in P$. But note that by $(1)$, for $a \in P$, $\Omega a = \{A \in \Omega: 1_{A}(a)=1\}$ which is clearly open in $\Omega$, as $\Omega$ has the subspace topology of $\{0,1\}^{G}$. This completes the proof. \hfill $\Box$

A consequence of Proposition \ref{discrete case} is that the semi-direct product groupoid $\Omega \rtimes P$ has a Haar system.
For $g \in G$, let $\delta_{g} \in C_{c}(\Omega \times P)$ be defined by $\delta_{g}(x,h)=1$ if $h=g$ and $\delta_{g}(x,h)=0$ if $h \neq g$. Then it is not difficult to show that there exists a representation $\pi:C_{c}(\Omega \times P) \to B(\clh)$ 
such that $\pi(\delta_{g})=W_{g^{-1}}$ for every $g \in G$. We will prove this in the topological case.

Now let us turn our attention to the topological case. Let $\chi$ be a character of the commutative $C^{*}$-algebra $\mathcal{A}$. The \textbf{support} of $\chi$, denoted $A_{\chi}$, is defined as follows: For $g \in G$, $g \notin A_{\chi}$ if and only if there exists an open set $U$ of $G$ containing $g$ such that $\chi(\int f(g)E_{g}dg)=0$ for every $f \in C_{c}(U)$. Here $C_{c}(U):=\{f \in C_{c}(G): supp(f) \subset U\}$. Note that $A_{\chi}$ is closed.

\begin{rmrk}
 \label{Support}
Let $\chi$ be a character of $\mathcal{A}$ and $A$ be its support. Then  for $g \in G$, $g \in A$ if and only if for every open set $U$ containing $g$, there exists $f \in C_{c}(U)$ such that $f \geq 0$ and $\chi(\int f(g)E_{g}dg) > 0$. 

\end{rmrk}

\begin{ppsn}
\label{supp}
 Let $\chi$ be a character of $\mathcal{A}$ and let $A$ be its support. Then 
\begin{enumerate}
 \item[(1)] $P^{-1} \subset A$ and $P^{-1}A \subset A$, 
 \item [(2)] the interior $Int(A)$ is dense in $A$, and
 \item [(3)] the boundary $\partial(A)$ has measure zero.
 \end{enumerate}
\end{ppsn}
\textit{Proof.} 
Let $a \in Int(P)$ and $U$ be an open set containing $a^{-1}$. Then $ U \cap Int(P)^{-1}$ is a non-empty open set containing $a^{-1}$. Choose $f \in C_{c}(G)$ such that $supp(f) \subset U \cap Int(P)^{-1}$, $f \geq 0$ and $\int f(g) dg=1$. Since $E_{g}=1$ for $g \in P^{-1}$, it follows that $\int f(g)E_{g}dg= \int f(g)dg=1$. Thus $\chi(\int f(g)E_{g}dg)=1$. This proves that $a^{-1} \in A$. As a consequence, $Int(P)^{-1} \subset A$. But $Int(P)^{-1}$ is dense in $P^{-1}$ and $A$ is closed. Hence $P^{-1} \subset A$.

For $f \in C_{c}(G)$ and $g \in G$, let $L_{g}(f) \in C_{c}(G)$ be defined by $L_{g}(f)(x)=f(g^{-1}x)$.

Let $g \in A$ be given and $a \in P$. Let $U$ be an open set containing $a^{-1}g$. Then $aU$ is open and contains $g$. Thus there exists $f \in C_{c}(aU)$ such that $f \geq 0$ and $\chi(\int f(g)E_{g}dg) > 0$. Let $\widetilde{f}=L_{a^{-1}}f$. Then 
 $\widetilde{f} \geq 0$ and $supp(\widetilde{f}) \subset U$.
Now 
\begin{align*}
     \int \widetilde{f}(g)E_{g}dg & = \int f(ag)E_{g}dg \\
                                  & = \int f(g)E_{a^{-1}g}dg \\
                                  & \geq \int f(g)E_{g}dg \textrm{~~(By Proposition \ref{commuting projections})}.
                                      \end{align*}
Hence $\chi(\int \widetilde{f}(g)E_{g}dg)  >0$. This implies that $a^{-1}g \in A$. Thus $P^{-1}A \subset A$. This proves $(1)$. Statements $(2)$ and $(3)$ follow immediately from Lemma 4.1 of \cite{Jean_Sundar}. This completes the proof. \hfill $\Box$

Before proceeding further, let us review the Vietoris topology. Let $X$ be a locally compact second countable Hausdorff space and let $d$ be a metric on $X$ inducing the topology. Let $\mathcal{C}(X)$ be the collection of closed subsets of $X$. Then $\mathcal{C}(X)$, endowed with the Vietoris topology, is compact and metrisable. We recall here  the convergence of sequences of elements in $\mathcal{C}(X)$.

Let $(A_{n})$ be a sequence of closed subsets of $X$. Define
 \begin{align*}
 \lim \inf A_{n} & = \{ x \in X: \lim \sup d(x,A_{n}) = 0 \}, \textrm{~and~} \\
  \lim\sup A_{n} & = \{x \in X: \lim \inf d(x,A_{n})=0 \}.
\end{align*}
Then $(A_{n})$ converges in $\mathcal{C}(X)$ if and only if $\lim \inf A_{n}=\lim \sup A_{n}$. If $\lim \inf A_{n}=\lim \sup A_{n}=A$, then $A_{n}$ converges to $A$. Observe that if $U \subset X$ is closed then the subset $\{A \in \mathcal{C}(X): A \cap U \neq \emptyset \}$ is open in $\mathcal{C}(X)$.

Consider $\mathcal{C}(G)$, the space of closed subsets of $G$, with the Vietoris topology. The group $G$ acts on $\mathcal{C}(G)$ on the right. For $A \in \mathcal{C}(G)$ and $g \in G$, define $Ag=\{ag:a \in A\}$. Let
\[
 \Omega_{u}:= \{A \in \mathcal{C}(G): P^{-1} \subset A \textrm{~and~}P^{-1}A \subset A\}.
\]
We leave it to the reader to verify that $\Omega_{u}$ is a closed, and hence a compact, subset of $\mathcal{C}(G)$. Clearly $\Omega_{u}$ is $P$-invariant.
The space $\Omega_{u}$ is first considered in \cite{Hilgert_Neeb}. 

\begin{ppsn}
\label{Open}
 The action $\Omega_{u} \times Int(P) \to \Omega_{u}$ is open.
\end{ppsn}
\textit{Proof.} Let $a \in Int(P)$.  It is enough to show that $\Omega_{u}Int(P)a$ is open in $\Omega_{u}$ (See Theorem 4.3, \cite{Jean_Sundar}). We claim that 
\[
 \Omega_{u}Int(P)a=\{A \in \Omega_{u}: A \cap Int(P)a \neq \emptyset \}
\]
which will imply that $\Omega_{u}Int(P)a$ is open. 

Let $A \in \Omega_{u}Int(P)a$. Then $A=Bba$ for some $B \in \Omega_{u}$ and $b \in Int(P)$. Since $e \in B$, it follows that $ba \in A$. Hence $A \cap Int(P)a$ is non-empty. Now suppose $A \in \Omega_{u}$ and $A \cap Int(P)a$ is non-empty. Choose $b \in Int(P)$ such that $ba \in A$. Since $P^{-1}A \subset A$, it follows that $P^{-1}ba \subset A$, equivalently $P^{-1} \subset Aa^{-1}b^{-1}$, and $P^{-1}Aa^{-1}b^{-1} \subset Aa^{-1}b^{-1}$. This proves that $B=Aa^{-1}b^{-1} \in \Omega_{u}$. Then $A=Bba \in \Omega_{u}Int(P)a$. This completes the proof. \hfill $\Box$

We summarise a few facts regarding the space $\Omega_{u}$ in the following remark. 
\begin{rmrk}
Note the following.
\label{Openness}
\begin{enumerate}
 \item [(1)] $\Omega_{u}  Int(P)a= \{ A \in \Omega_{u}: a \in Int(A) \}$. If $A \cap Int(P)a \neq \emptyset$ then $a \in Int(P)^{-1}A$ which is open and contained in $A$. Thus $a \in Int(A)$. Now suppose $a \in Int(A)$ then $Int(A) \cap Pa$ is non-empty. Since $Int(P)a$ is dense in $Pa$, it follows that $Int(A) \cap Int(P)a$ is non-empty and hence $A \cap Int(P)a$ is non-empty.
\item [(2)] If $A \in \Omega_{u}$ then $Int(A)$ is dense in $A$ and the boundary $\partial(A)$ has measure zero. This follows from Lemma 4.1 of \cite{Jean_Sundar}
 \item [(3)] Let $A \in \Omega_{u}$ and $g \in G$. Then $(A,g) \in \Omega_{u} \rtimes P$ if and only if $Ag \in \Omega_{u}$ if and only if $g^{-1} \in A$. We leave this verification to the reader.
             
\item [(4)] The map $\Omega_{u} \ni A \to 1_{A} \in L^{\infty}(G)$ is continuous and injective and hence an embedding. Here $L^{\infty}(G)$ is given the weak $^{*}$-topology. Let $\mathcal{G}_{u}:=\Omega_{u} \rtimes P$. Then Proposition \ref{Open} implies that $\mathcal{G}_{u}$ has a Haar system. Moreover a Haar system on $\mathcal{G}_{u}$ is given by $(1_{Q_{A}}(g)dg)_{A \in \Omega_{u}}$. For $A \in \Omega_{u}$, observe that $Q_{A}:=\{g \in G: (A,g) \in \mathcal{G}_{u}\}$ is $A^{-1}$. 
 
 By the definition of a Haar system, it follows that for $f \in C_{c}(\mathcal{G}_{u})$,  $\Omega_{u} \ni A \to \int f(x,g)1_{A}(g^{-1})dg$ is continuous. In particular, for $f \in C_{c}(G)$, the function $\Omega_{u} \ni A \to \int f(g)1_{A}(g^{-1})dg$ is continuous. As a consequence, the map $\Omega_{u} \ni A \to 1_{A} \in L^{\infty}(G)$ is continuous.

Suppose $A,B \in \Omega_{u}$ such that $1_{A}=1_{B}$ in $L^{\infty}(G)$. Then $A\backslash B$ and $B\backslash A$ has measure zero. If $A \backslash  B$ is non-emtpy then $Int(A) \backslash B$ is non-empty since $Int(A)$ is dense in $A$. But $Int(A) \backslash B$ is open and hence cannot have measure zero. Thus $A\backslash B=\emptyset$. Similarly $B\backslash A = \emptyset$. Hence $A=B$. This proves the map $\Omega_{u} \ni A \to 1_{A} \in L^{\infty}(G)$ is injective. Thus we can consider $\Omega_{u}$ as a compact subset of $L^{\infty}(G)$.
\end{enumerate}
\end{rmrk}

Let $V:P \to B(\clh)$ be an isometric representation with commuting range projections.  Denote the  commutative $C^{*}$-algebra generated by $\{ \int f(g)E_{g}dg: f \in L^{1}(G)\}$ by $\cla$  and let  $\Omega$ be the spectrum of $\mathcal{A}$. For $f \in L^{1}(G)$, let $E_{f}:=\int f(g)E_{g} dg$.  

\begin{ppsn}
\label{dfn of spectrum}
 Let $\chi$ be a character of $\mathcal{A}$ and let $A$ be its support. Let $f \in C_{c}(G)$. Then
\begin{enumerate}
 \item[(1)]  $\chi \Big( \int f(g)1_{A^{c}}(g)E_{g}dg \Big)=0$.
 \item [(2)] if $supp(f) \subset Int(A)$ then $\chi\Big(\int f(g)E_{g}dg\Big)=\int f(g)dg$, and
 \item [(3)] we have the equality $\chi \Big( \int f(g)E_{g} dg \Big)= \int f(g)1_{A}(g)dg$.
\end{enumerate}
\end{ppsn}
\textit{Proof.} First observe that if $supp(f) \subset A^{c}$, where $A^{c}$ denotes the complement of $A$, then $\chi(E_{f})=0$. This follows from the definition of $A$ and by a partition of unity argument. 

Now write $A^{c}=\displaystyle \bigcup_{n=1}^{\infty}{K_{n}}$ with $K_{n}$ compact and $K_{n}$ increasing. This is possible as $A^{c}$ is open. Choose $\phi_{n} \in C_{c}(G)$ such that $0 \leq \phi_{n} \leq 1$, $\phi_{n}=1$ on $K_{n}$ and $supp(\phi_{n}) \subset A^{c}$. Note that  $\phi_{n} \to 1_{A^{c}}$ pointwise. Hence $f\phi_{n}$ converges to $f1_{A^{c}}$ in $L^{1}(G)$. This implies that $E_{f\phi_n}$ converges to $E_{f1_{A^{c}}}$. Since $\chi(E_{f\phi_{n}})=0$, it follows that $\chi(E_{f1_{A^{c}}})=0$. This proves $(1)$.

Let $f \in C_{c}(G)$ be such that $supp(f) \subset Int(A)$. Let $g \in supp(f)$. Then $Int(A) \cap Pg$ is non-empty. Since $Int(P)$ is dense in $P$, it follows that $Int(A)\cap Int(P)g$ is non-empty. Let $s \in Int(P)$ be such that $sg \in Int(A)$.
Then $(sg)g^{-1} \in Int(P)$. Since $Int(P)$ is open, we can choose open sets $U$ and $V$ contained in $Int(A)$, with compact closures, such that $(g,sg) \in U \times V \subset Int(A) \times Int(A)$ and $VU^{-1} \subset Int(P)$. Then by Proposition \ref{commuting projections}, it follows that for $g_{1} \in V$ and $g_{2} \in U$, $E_{g_{1}}E_{g_{2}}=E_{g_{1}}$.

Since $supp(f)$ is compact, it follows that there exists finitely many non-empty open sets $(U_{i})_{i=1}^{n}$ and $(V_{i})_{i=1}^{n}$ with compact closures, contained in $Int(A)$, such that $supp(f) \subset \bigcup_{i=1}^{n}U_{i}$ and $V_{i}U_{i}^{-1} \subset Int(P)$. A partition of unity argument allows us to write $f$ as $f=\sum_{i=1}^{n}f_{i}$ with $supp(f_{i}) \subset U_{i}$. Thus to prove $(2)$, it is enough to show  $\chi(E_{f_{i}})=\int f_{i}(g)dg$.

Since $V_{i}$ is a non-empty open set contained in $A$, by Remark \ref{Support}, it follows that there exists $\phi_{i} \in C_{c}(G)$ such that $supp(\phi_{i}) \subset V_{i}$ and $\chi(E_{\phi_{i}}) \neq 0$. Observe the following
\begin{align*}
 E_{\phi_{i}}E_{f_{i}}& = \int_{V_{i}\times U_{i}} \phi_{i}(g_{1})f_{i}(g_{2})E_{g_{1}}E_{g_{2}}dg_{1}dg_{2} \\
                      & =\int_{V_{i}\times U_{i}} \phi_{i}(g_{1})f_{i}(g_{2})E_{g_{1}}dg_{1}dg_{2} \\
                      & = \Big ( \int f_{i}(g_{2})dg_{2} \Big)\int \phi_{i}(g_{1})E_{g_{1}}dg_{1} \\
                      & = \Big ( \int f_{i}(g_{2})dg_{2} \Big)E_{\phi_{i}}.
\end{align*}
Since $\chi$ is multiplicative, it follows that $\chi(E_{\phi_{i}})\chi(E_{f_{i}})= \Big(\int f_{i}(g)dg\Big)\chi(E_{\phi_{i}})$. Now $\chi(E_{\phi_{i}}) \neq 0$. Hence $\chi(E_{f_{i}})=\int f_{i}(g)dg$. This proves $(2)$.

Now let $f \in C_{c}(G)$ be given. By $(1)$, it follows that $\chi(E_{f})=\chi(E_{f1_{A}})$. But since the boundary of $A$ has measure zero, it follows that $1_{Int(A)}=1_{A}$ a.e. Thus $\chi(E_{f})=\chi(E_{f1_{Int(A)}})$. Write $Int(A)=\bigcup_{n}K_{n}$ with $K_{n}$ compact and $K_{n}$ increasing. Choose $\phi_{n} \in C_{c}(G)$ such that $\phi_{n}=1$ on $K_{n}$ and $supp(\phi_{n}) \subset Int(A)$. Then $\phi_{n} \to 1_{Int(A)}$ pointwise and hence $f\phi_{n}$ converges to $f1_{IntA}$ in $L^{1}(G)$. Note that $supp(f\phi_{n}) \subset Int(A)$. Now calculate, as follows, to find that
\begin{align*}
\chi(E_{f}) & = \chi(E_{f 1_{Int(A)}})\\
            & = \lim_{n} \chi (E_{f\phi_{n}})\\
            & = \lim_{n} \int f(g)\phi_{n}(g)dg ~~(\textrm{by~} (2))\\
            &  =  \int f(g)1_{Int(A)}(g)dg \\
            & = \int f(g)1_{A}(g)dg~~(\textrm{~Since $1_{A}=1_{Int(A)}$ in $L^{\infty}(G)$).} 
 \end{align*}
This proves $(3)$. This completes the proof. \hfill $\Box$.

\begin{ppsn}
 \label{Spectrum}
For $\chi \in \Omega$, let $A_{\chi}$ be its support. Then the map $\Omega \ni \chi \to A_{\chi} \in \Omega_{u}$ is one-one, continuous and $P$-equivariant. Consequently, the action of $P$ on $\Omega$ is injective.
\end{ppsn}
\textit{Proof.} By Proposition \ref{supp}, it follows that $A_{\chi} \in \Omega_{u}$ if $\chi \in \Omega$. For $f \in C_{c}(G)$, by Proposition \ref{dfn of spectrum}, $\chi(E_{f})=\int f(g)1_{A_{\chi}}(g)dg$. Hence the map $\Omega \ni \chi \to A_{\chi} \in L^{\infty}(G)$ is one-one and continuous where $L^{\infty}(G)$ is given the weak $^{*}$-topology. By part $(4)$ of Remark \ref{Openness}, it follows that $\Omega \ni \chi \to A_{\chi} \in \Omega_{u}$ is one-one and continuous.

Let $f \in C_{c}(G)$, $a \in P$ and $\chi \in \Omega$ and $A$ be the support of $\chi$. Observe that
\begin{align*}
 (\chi.a)(\int f(g)E_{g}dg)&=\chi(\int f(g)V_{a}^{*}E_{g}V_{a})dg \\
                           & = \chi (\int f(g)E_{ga^{-1}}dg) \\
                           & = \chi (\int f(ga)\Delta(a)E_{g}dg) \\
                           & = \int f(ga)\Delta(a)1_{A}(g)dg \\
                           &= \int f(g)1_{A}(ga^{-1})dg \\
                           & = \int f(g)1_{Aa}(g) dg.
\end{align*}
Hence the support of $\chi.a$ is $Aa$. Thus the map $\Omega \ni \chi \to A_{\chi} \in \Omega_{u}$ is a continuous $P$-equivariant embedding. This completes the proof. \hfill $\Box$.

Thus we can and will consider $\Omega$ as a subset of $\Omega_{u}$ with the subspace topology.

\section{Haar system on $\Omega \rtimes P$}

In this section, we show that the semi-direct product $\Omega \rtimes P$ admits a Haar system. We prove that the action $\Omega \times Int(P) \to \Omega$ is open. To prove this, we need an analogue of Proposition \ref{discrete case} in the topological setting.

\begin{ppsn}
\label{ouvert}
 Let $a \in Int(P)$ and $A \in \Omega$. Then $a \in A$ if and only if $Aa^{-1} \in \Omega$.
\end{ppsn}
\textit{Proof.} Let $a \in Int(P)$ and $A \in \Omega$ be given. Suppose $B:=Aa^{-1} \in \Omega$.  Since $e \in B$, it follows that $a \in A$. Now suppose $a \in A$. In addition, assume  that $a \in Int(A)$.  Let $\chi$ be the character defining $A$. Then for $f \in C_{c}(G)$, \[\chi(\int f(g)E_{g}dg)=\int f(g)1_{A}(g)dg.\] Choose a decreasing sequence of open sets $(U_{n})$ in $G$ such that
\begin{enumerate}
 \item [(1)] the intersection $\bigcap_{n=1}^{\infty}U_{n}=\{a\}$,
 \item [(2)] if $U$ is open in $G$ and $a \in U$ then there exists $N$ such that $U_{n} \subset U$ for $n \geq N$, and
 \item [(3)] for every $n$, $U_{n} \subset Int(P) \cap Int(A)$.
\end{enumerate}
This is possible, for we can choose a metric and let $(U_{n})$  be the open balls containing $a$ with $diam(U_{n})\to 0$. For every $n \in \mathbb{N}$, choose $f_{n} \in C_{c}(G)$ such that $ f_{n} \geq 0$, $\int f_{n}(g)dg=1$ and $supp(f_{n}) \subset U_{n}$. Note that $\chi(E_{f_{n}})=\int f_{n}(g)1_{A}(g)dg = 1$ since $supp(f_{n}) \subset Int(A)$.

Let $\phi_{n}$ be the linear functional on the commutative $C^{*}$-algebra $\mathcal{A}$ defined by \[
                                                                           \phi_{n}(T)=\chi \big ( \int_{b \in P} f_{n}(b)V_{b}TV_{b}^{*}db \big).
                                                                          \]
Note that $\phi_{n}$ is well defined by Lemma \ref{Transfer operator} and is clearly positive. Note that 
\begin{align*}
 \phi_{n}(1)&= \chi \Big(\int_{b \in supp(f_{n})}f_{n}(b)E_{b}db \Big) \\
            & = \int_{b \in supp(f_{n})}f_{n}(b)1_{A}(b)db \\
            & = \int_{b \in supp(f_{n})}f_{n}(b)db ~ (\textrm{ Since $supp(f_{n}) \subset Int(A)$}) \\
            & =1.
\end{align*}
Thus $\phi_{n}$ is a state for every $n$.  But the set of states on a unital $C^{*}$-algebra is weak$^{*}$-compact. By choosing a subsequence if necessary we can assume without loss of generality that $(\phi_{n})$ converges in the weak$^{*}$-topology and let $\phi$ be its limit. 

Recall that for $a \in P$, $\alpha_{a}:\mathcal{A} \to \mathcal{A}$ is given by $\alpha_{a}(T)=V_{a}^{*}TV_{a}$. By Proposition \ref{Spectrum}, it follows that for every $a \in \mathcal{A}$, $\alpha_{a}$ is surjective.

\textit{Claim:}  $\phi \circ \alpha_{a}=\chi$. It is enough to show that $(\phi \circ \alpha_{a})(T)=\chi(T)$ for $T \in \widetilde{\mathcal{A}}$. Let $T=E_{\psi}$ for some $\phi \in C_{c}(G_{m})$.

Observe that for $n \in \mathbb{N}$, 
\begin{align*}
 &\phi_{n}(\alpha_{a}(E_{\psi})) \\
&= \chi(\int_{b \in P}f_{n}(b)V_{b}\alpha_{a}(E_{\psi})V_{b}^{*}db ) \\
                               & = \chi (\int_{b \in P}f_{n}(b)\psi(g_{1},g_{2},\cdots,g_{m})V_{b}V_{a}^{*}E_{g_{1}}\cdots E_{g_{m}}V_{a}V_{b}^{*} db dg_{1}\cdots dg_{m} ) \\
                               & = \chi ( \int_{b \in P}f_{n}(b)\psi(g_{1},g_{2},\cdots,g_{m})E_{b}E_{g_{1}a^{-1}b}\cdots E_{g_{m}a^{-1}b} db dg_{1} \cdots dg_{m} ) \\
                               & = \chi (\int_{b \in P}\Delta(b^{-1}a)^{m}f_{n}(b)\psi(g_{1}b^{-1}a,g_{2}b^{-1}a,\cdots,g_{m}b^{-1}a)E_{b}E_{g_{1}}E_{g_{2}}\cdots E_{g_{m}} db dg_{1} dg_{2} \cdots dg_{m}).\\
\end{align*}
Let $\epsilon >0$ be given. Since $\psi$ is continuous and compactly supported, it follows that there exists an open set $U$ such that $a \in U$ and for $b \in U$ and $g_{1},g_{2},\cdots,g_{m} \in G$,
\[
 \int |\Delta(b^{-1}a)^{m}\psi(g_{1}b^{-1}a,g_{2}b^{-1}a,\cdots,g_{m}b^{-1}a)-\psi(g_{1},g_{2},\cdots,g_{m})| dg_{1}dg_{2}\cdots dg_{m} \leq \epsilon.
\]

Choose $N \geq 1$ such that for $n \geq N$, $U_{n} \subset U$. Then for $n \geq N$, $supp(f_{n}) \subset U$. 
Note that for $n \geq N$, 
\begin{align*}
& |\phi_{n}(\alpha_{a}(E_{\psi}))-\chi(E_{\psi})| \\
&= | \phi_{n}(\alpha_{a}(E_{\psi}))-\chi(E_{f_{n}}E_{\psi})|\\
& \leq \int_{b \in U_{n}}f_{n}(b)\Big(\int |\Delta(b^{-1}a)^{m}\psi(g_{1}b^{-1}a,\cdots,g_{m}b^{-1}a)-\psi(g_{1},g_{2},\cdots,g_{m})| dg_{1} \cdots dg_{m} \Big)db \\ 
& \leq \epsilon \int _{b \in U_{n}}f_{n}(b) \\
& \leq \epsilon.
\end{align*}
 Thus it follows that $\phi_{n}(\alpha_{a}(E_{\psi}))\to \chi(E_{\psi})$ and hence $\phi \circ \alpha_{a} = \chi$. This proves the claim.

Since $\alpha_{a}$ is surjective on $\mathcal{A}$, it follows that $\phi$ is a character of $\mathcal{A}$. Let $B \in \Omega$ be the support of $\phi$. Then $\phi \circ \alpha_{a}=\chi$ translates to the equation $Ba=A$. Thus $Aa^{-1} \in \Omega$.
Now suppose $a \in A$. Let $(s_{n})$ be a sequence in $Int(P)$ converging to the identity element $e$. Then $s_{n}^{-1}a \in Int(P)$ eventually, for $s_{n}^{-1}a \to a$ and $a \in Int(P)$.   But $Int(P)^{-1}A \subset Int(A)$. Hence $s_{n}^{-1}a \in Int(A)$. By what we have proved, it follows that $Aa^{-1}s_{n} \in \Omega$ eventually. However $\Omega$ is a compact subset of $\Omega_{u}$ and $(Aa^{-1}s_{n})$ converges to $Aa^{-1}$. From this we conclude that $Aa^{-1} \in \Omega$. This completes the proof. \hfill $\Box$

Just like in the discrete case, we have the following theorem.

\begin{ppsn}
\label{Haar system}
 Let $A \in \Omega$ and $g \in G$. Then $Ag \in \Omega$ if and only if $g^{-1} \in A$. Also the semi-direct product groupoid $\Omega \rtimes P$ has a Haar system.
\end{ppsn}
\textit{Proof.} Let $g \in G$ and $A \in \Omega$ be given. Suppose $Ag \in \Omega$. Since $e \in Ag$, it follows that $g^{-1} \in A$.
Now suppose $g^{-1} \in A$. As $G=(Int(P))(Int(P))^{-1}$, write $g=ab^{-1}$ with $a,b \in Int(P)$. Then $ba^{-1} \in A$ or $b \in Aa \in \Omega$. By Proposition \ref{ouvert}, it follows that $Aab^{-1} \in \Omega$. Hence $Ag \in \Omega$. 

To prove that $\Omega \rtimes P$ has a Haar system, it is enough to show that the action $\Omega \times P \to \Omega$ is open.
By Theorem 4.3 of \cite{Jean_Sundar}, it is enough to show that $\Omega Int(P)a$ is open in $\Omega$ for every $a \in P$. Let $a \in P$ be given.

 \textit{Claim:} $\Omega Int(P)a=\{A \in \Omega: A \cap Int(P)a \neq \emptyset\}$.

 Suppose $A \cap Int(P)a$ is non-empty. Then there exists $s \in Int(P)$ such that $sa \in A$. By Proposition \ref{ouvert}, $B:=Aa^{-1}s^{-1} \in \Omega$. Thus $A=Bsa \in \Omega Int(P)a$. Suppose $A \in \Omega Int(P)a$. Then $A=Bsa$ for some $B \in \Omega$ and $s \in Int(P)$. Since $e \in B$, it follows that $sa \in A$. Thus $A \cap Int(P)a$ is non empty. This proves the claim.

The set $\{A \in \mathcal{C}(G): A \cap Int(P)a \neq \emptyset \}$ is open in $\mathcal{C}(G)$, when $\mathcal{C}(G)$ is given the Vietoris topology. This implies that $\Omega Int(P)a$ is open in $\Omega$. This completes the proof. \hfill $\Box$




\begin{rmrk}
\label{restriction}
Consider the groupoid $\Omega_{u} \rtimes P$. Then by Proposition \ref{Haar system} and Statement $(3)$ of Remark \ref{Openness}, it follows that $\Omega$ is an invariant subset of $\Omega_{u}$. Moreover the groupoid $\Omega \rtimes P$ is just the restriction $\Omega_{u} \rtimes P|_{\Omega}$.
\end{rmrk}

We end this section by describing $\Omega$ in the case of the Wiener-Hopf representation.  Recall that the Wiener-Hopf represention $V:P \to B(L^{2}(P))$ is given by the formula: For $a \in P$ and $\xi \in L^{2}(P)$, \begin{displaymath}
\begin{array}{lll}
V_{a}(\xi)(x)&:=& \left\{\begin{array}{lll}
                               \xi(xa^{-1})\Delta(a)^{\frac{-1}{2}} &~if~x \in Pa \\
                                 0 & ~if~x \notin Pa.
                                 \end{array} \right. 
\end{array}
\end{displaymath}
Here $\Delta$ denotes the modular function of the group $G$. Note that for $g \in G$ and $\xi \in L^{2}(P)$, $W_{g}$ is given by \begin{displaymath}
\begin{array}{lll}
W_{g}(\xi)(x)&:=& \left\{\begin{array}{lll}
                               \xi(xg^{-1})\Delta(g)^{\frac{-1}{2}} &~if~xg^{-1} \in P \\
                                 0 & ~if~xg^{-1} \notin P.
                                 \end{array} \right. 
\end{array}
\end{displaymath}
Let $M:L^{\infty}(P) \to B(L^{2}(P))$ be the multiplication representation. Observe that for $g \in G$, $E_{g}=W_{g}W_{g}^{*}=M(1_{P.g})$ where $P.g:= \{xg:x \in P\} \cap P$. Denote the algebra of bounded continuous functions on $P$ by $C_{b}(P)$. Since $\overline{Int(P)}=P$, it follows that $M$ is a  faithful representation of  $C_{b}(P)$. For $f \in C_{c}(G)$, let $1_{P}*f \in C_{b}(P)$ be defined by \begin{align*}
                                     1_{P}*f(t)&=\int 1_{P}(ts)f(s^{-1})ds\\
                                               &= \int 1_{P}(ts^{-1})f(s)\Delta(s)^{-1}ds\\
                                               &=\int 1_{P^{-1}t}(s)f(s)\Delta(s)^{-1}ds
                                    \end{align*}
 Observe that given $a \in P$, there exists $f \in C_{c}(G)$ such that $(1_{P}*f)(a)=1$. For, if $a \in P$, choose $f \in C_{c}(G)$ such that $supp(f) \subset Int(P)^{-1}a$ and $\int f(s)\Delta(s)^{-1}ds=1$. For such an $f$, $(1_{P}*f)(a)=1$. 
                                    
Now let $f \in C_{c}(G)$ and $\xi \in L^{2}(P)$ be given. Calculate as below to find that
\begin{align*}
 <(\int f(g)E_{g}dg)\xi,\xi > & = \int f(g)<E_{g}\xi,\xi>dg \\
                              & = \int f(g) \Big( \int_{x \in P} 1_{Pg}(x)|\xi(x)|^{2}dx \Big) dg \\
                              & = \int_{x \in P} \Big( \int f(g)1_{Pg}(x) dg\Big) |\xi(x)|^{2} dx \\
                              & = \int_{x \in P} \Big(\int f(g)1_{P}(xg^{-1})dg\Big) |\xi(x)|^{2} dx \\
                              & = \int_{x \in P} \Big(\int 1_{P}(xg)f(g^{-1})\Delta(g)^{-1}dg\Big) |\xi(x)|^{2} dx \\
                              & = \int_{x \in P} (1_{P}*\hat{f})(x)|\xi(x)|^{2}dx \\
                                 & = <  M(1_{P}*\hat{f})\xi,\xi> .                   
\end{align*}
where $\hat{f}(g)=f(g)\Delta(g)$. Thus  the $C^{*}$-algebra generated by $\{\int f(g)E_{g}dg: f \in C_{c}(G)\}$, is isomorphic to the $C^{*}$- subalgebra of $C_{b}(P)$ generated by $\{1_{P}*f: f \in C_{c}(G)\}$. 
Thus for $a \in P$, there exists a character $\chi_{a}$ of $\mathcal{A}$ such that
\begin{align*}
\chi_{a}(\int f(g)E_{g}dg)&= (1_{P}*\hat{f})(a) \\
                                       & = \int 1_{P^{-1}a}(s)\hat{f}(s)\Delta(s)^{-1} ds \\
                                        &= \int 1_{P^{-1}a}(s)f(s) ds
\end{align*}
This implies that the support of $\chi_{a}$ is $P^{-1}a$.  Also $\{\chi_{a}: a \in P\}$ separates the elements of $C_{b}(P)$ and hence those of $\mathcal{A}$. This implies that $\{P^{-1}a: a \in P\}$ is dense in $\Omega$. As a consequence, it follows that $\Omega$ is the closure of $\{P^{-1}a: a \in P\}$ in the space of closed subsets of $G$ w.r.t. the Vietoris topology. The 'compactification' $\Omega$ of $P$ is called the Wiener-Hopf compactification and is considered in \cite{Renault_Muhly} and in \cite{Jean_Sundar}.

\section{Covariant representations}
In this section, let $X$ be a compact Hausdorff space and assume that $P$ acts on $X$ on the right injectively.  Let $X_{0}:=X Int(P)$. We also assume that the semi-direct product $\mathcal{G}:=X \rtimes P$ admits a Haar system.  Let $Y$ be a dilation of $X$, as explained in Section 1, on which the group $G$ acts.  For $y \in Y$, let $Q_{y}:=\{g \in G: y.g \in X\}$. Recall that for $y \in Y$ and $g \in G$, $1_{Q_{y}}(g)=1_{X}(yg)$ and $1_{Int(Q_{y})}(g)=1_{X_{0}}(yg)$. Also note that for every $y \in Y$, $Q_{y}$ is closed and $Q_{y}P \subset Q_{y}$. Thus by Lemma 4.1 of \cite{Jean_Sundar}, it follows that  for every $y \in Y$, the boundary of $Q_{y}$ has measure zero and $\overline{Int(Q_{y})}=Q_y$.

To state the next lemma, we need to fix some notations. Let $a \in P$ and let $(U_{n})$ be a decreasing sequence of open subsets of $G$ such that $\displaystyle \bigcap_{n=1}^{\infty}U_{n}=\{a\}$ and if $U$ is open and contains $a$ then  $U_{n} \subset U$ eventually. Note that for every $n$, $U_{n} \cap Pa$ is non-empty. Hence $U_{n} \cap Int(P)a$ is non-empty for every $n$.  Choose $f_{n} \in C_{c}(G)$ such that $ f_{n} \geq 0$, $  \int f_{n}(g)dg=1$ and $supp(f_{n}) \subset U_{n} \cap Int(P)a$. For $x \in X$, let \[F_{n}(x)=\int f_{n}(g)1_{X}(xg^{-1})dg.\] Then $F_{n} \in C(X)$. The continuity of $F_n$ follows from the fact that $(1_{Q_{x}}(g)dg)_{x \in X}$ is a Haar system on $X \rtimes P$. Observe that $(F_{n})$ is uniformly bounded.

\begin{lmma}
\label{crying out lemma}
The sequence $F_{n}$ converges pointwise to $1_{X_{0}a}$.
\end{lmma}
\textit{Proof.} Since for $x \in X$, $1_{Q_{x}}=1_{Int(Q_{x})}$ a.e., it follows that $F_{n}$ is given by the equation
\[F_{n}(x)=\int f_{n}(g)1_{X_{0}}(xg^{-1})\]
 for $x \in X$.  
Let $x \in X$. From $Q_{x}P \subset Q_{x}$, it is easily verifiable that $a^{-1} \in Int(Q_{x})$ if and only if $Int(Q_{x}) \cap a^{-1}P^{-1}$ is non-empty. Suppose $a^{-1} \in Int(Q_{x})$ i.e. $xa^{-1} \in X_{0}$. 

 Let $U:=\{g \in G: xg^{-1} \in X_{0}\}$. Then $U$ is open and contains $a$. Thus there exists $N$ such that $n \geq N$ implies $supp(f_{n}) \subset U$ eventually. Then for $n \geq N$, $F_{n}(x)=\int f_{n}(g)dg =1$. 
Now suppose $xa^{-1} \notin X_{0}$ i.e. $a^{-1} \notin Int(Q_{x})$.  Then $Int(Q_{x})^{-1} \cap Pa$ is empty. Thus for $g \in Pa$, $xg^{-1} \notin X_{0}$. Since $supp(f_{n}) \subset Pa$, it follows that $F_{n}(x)=0$. This proves that  $(F_{n})$ converges pointwise to $1_{X_{0}a}$. This completes the proof. \hfill $\Box$

\begin{lmma}
\label{approximate identity}
 There exists a sequence $(s_{n})$ in $Int(P)$ such that $s_{n+1}^{-1}s_{n} \in Int(P)$ and $(s_{n})$ converges to the identity element $e$
\end{lmma}
\textit{Proof.} Let $(U_{n})$ be a countable base (of open sets) at $e$. We can assume that $U_{n}$ is decreasing. Now $U_{1} \cap P$ contains $e$ and is non-empty. Since $Int(P)$ is dense in $P$, it follows that $U_{1} \cap Int(P)$ is non-empty. Pick $s_{1} \in U_{1} \cap Int(P)$. Now suppose that $s_{1}, s_{2}, \cdots, s_{n}$ are chosen such that $s_{k} \in Int(P) \cap U_{k}$ for $1 \leq k \leq n$ and $s_{k+1}^{-1}s_{k} \in Int(P)$ for $1 \leq k \leq n-1$. Since $e \in s_{n}(Int(P))^{-1} \cap U_{n+1}$, it follows that $s_{n}Int(P)^{-1} \cap U_{n+1} \cap P$ is non-empty. But $\overline{Int(P)}=P$. Thus $s_{n}Int(P)^{-1} \cap U_{n+1} \cap Int(P)$ is non-empty. Let $s_{n+1} \in s_{n}Int(P)^{-1} \cap U_{n+1} \cap Int(P)$.

Then it is clear that the sequence $(s_{n})$ constructed as above converges to $e$ and $s_{n+1}^{-1}s_{n} \in Int(P)$ for every $n$. This completes the proof. \hfill $\Box$.
 
Consider a sequence  $(s_{n})$  as in Lemma \ref{approximate identity} converging to the identiy $e$. Let $a \in Int(P)$ and set $t_{n}:=s_{n}^{-1}a$. Then observe that $t_{n+1}t_{n}^{-1} \in Int(P)$. Since $a \in Int(P)$, $Int(P)$ is open and $(t_{n})$ converges to $a$, we can assume  without loss of generality that $t_{n} \in Int(P)$ for every $n$. With this notation, we have the following lemma.

\begin{lmma}
\label{decreasing sequence}
The sequence $(1_{X_{0}t_{n}})$ decreases pointwise to $1_{Xa}$.
\end{lmma}
\textit{Proof.} Since $t_{n+1}t_{n}^{-1} \in Int(P)$, it follows that $X_{0}t_{n+1} \subset X_{0}t_{n}$ for every $n$.  Since $at_{n}^{-1} \in Int(P)$, it follows that $Xa \subset X_{0}t_{n}$ for every $n$. Thus $\displaystyle Xa \subset \bigcap_{n=1}^{\infty}X_{0}t_{n}$.
Now suppose $y \in X_{0}t_{n}$ for every $n$. Then $yt_{n}^{-1} \subset X_{0}$ for every $n$. Note that $(yt_{n}^{-1}) \to ya^{-1}$. Since the closure of $X_{0}$ in $Y$ is $X$, it follows that $ya^{-1} \in X$. Hence $y \in Xa$. This proves that $\displaystyle Xa=\bigcap_{n=1}^{\infty}X_{0}t_{n}$. This completes the proof. \hfill $\Box$

Let $B(X)$ be the space of bounded Borel measurable functions on $X$. For $\phi$ in $B(X)$ and $g \in G$, let $
R_{g}(\phi)$  be defined by 
\begin{displaymath}
\begin{array}{lll}
R_{g}(\phi)(x)&:=\left\{\begin{array}{lll}
         \phi(x.g) & \textrm{if} & x.g \in X \\
                                 0 & \textrm{if} & x.g \notin X 
                                 \end{array} \right. 
\end{array}
\end{displaymath}
Then $R_{g}(\phi) \in B(X)$.

\begin{dfn}
Let $\pi:C(X) \to B(\clh)$ be a unital $*$-representation and $V:P \to B(\clh)$ be an isometric representation with commuting range projections. Denote the extension of $\pi$ to $B(X)$, obtained via the Riesz representation theory, by $\pi$ itself [See \cite{Arv}]. For $g=ab^{-1}$, let $W_{g}:=V_{b}^{*}V_{a}$.  The pair $(\pi,V)$ is said to be a covariant representation of $(X,P)$ if for $\phi \in B(X)$, \[W_{g}\pi(\phi)W_{g}^{*}=\pi(R_{g^{-1}}(\phi)).\]
\end{dfn}

\begin{rmrk}
\label{covariance}
Since $G=(Int(P))(Int(P))^{-1}$, it follows that $(\pi,V)$ is a covariant representation if and only if $V_{a}^{*}\pi(\phi)V_{a}=\pi(R_{a}(\phi))$ and $V_{a}\pi(\phi)V_{a}^{*}=\pi(R_{a^{-1}}(\phi))$ for $\phi \in B(X)$ and $a\in Int(P)$. We leave this verification to the reader. 
\end{rmrk}

We fix a few notations that will be useful for the rest of this section.

\textbf{Notations:} Let $Y$ be the dilation of $X$, as explained in Section 1, on which $G$ acts. Then $\mathcal{G}:=X \rtimes P$ is a closed subset of $X \times G$ and also of $Y \rtimes G$. For $\phi \in C_{c}(Y)$, we let $\widehat{\phi} \in C(X)$ be the restriction.  Define  $\pi_{Y}:C_{c}(Y) \to B(\clh)$ by $\pi_{Y}(\phi)=\pi(\widehat{\phi})$.
For $\phi \in C_{c}(Y)$ and $g \in G$, let $R_{g}(\phi) \in C_{c}(Y)$ be given by $R_{g}(\phi)(y)=\phi(y.g)$ for $y \in Y$. 

For $\psi \in C_{c}(Y \rtimes G)$ and $g \in G$, let $\psi_{g} \in C_{c}(Y)$ be defined by $\psi_{g}(y)=\psi(y,g)$. For $\chi \in C_{c}(\mathcal{G})$ and $g \in G$, let $\chi_{g} \in B(X)$ be defined by $\chi_{g}(x)=\chi(x,g)$ if $(x,g) \in \mathcal{G}$ and $\chi_{g}(x)=0$ if $(x,g) \notin \mathcal{G}$. 

Let $\pi:C(X) \to B(\clh)$ be a unital $*$-representation. For $\xi \in \clh$, let $d\mu_{\xi,\xi}$ be the probability measure on $X$ such that \[\int \phi(x)d\mu_{\xi,\xi}(x)= < \pi(\phi)\xi,\xi> ~ for ~ \phi \in  C(X).\] The same equality holds for $\phi \in B(X)$.

\begin{ppsn}
\label{covariant}
 Let $\pi:C(X) \to B(\clh)$ be a unital $*$-representation and $V:P \to B(\clh)$ be an isometric representation with commuting range projections.
Then the following are equivalent.
\begin{enumerate}
 \item [(1)] The pair $(\pi,V)$ is a covariant representation.
 \item [(2)] For $a \in Int(P)$ and $\phi \in C(X)$, $V_{a}^{*}\pi(\phi)V_{a}=\pi(R_{a}(\phi))$ and $\pi(1_{X_{0}a})=E_{a}$.
\end{enumerate}
\end{ppsn}
\textit{Proof.} For $a \in P$, let $\sigma_{a}:X \to X$ be the map sending $x \to xa$. 

Suppose $(\pi,V)$ is a covariant representation. Then the covariance relation implies that for $g \in G$, $E_{g}:=W_{g}W_{g}^{*}=\pi(1_{Xg \cap X})$.

Let $f \in C_{c}(G)$. Set $F(x)=\int f(g)1_{X}(xg^{-1})dg$. Then $F \in C(X)$. Now calculate to find that 
\begin{align*}
< \pi(F)\xi,\xi> &= \int F(x) d\mu_{\xi,\xi} (x)\\
                           & = \int f(g)1_{Xg}(x)dg d \mu_{\xi,\xi}(x) \\
                            & = \int f(g) \Big (\int 1_{Xg}(x)d\mu_{\xi,\xi}(x) \Big) dg \\
                             & = \int f(g) < \pi(1_{Xg\cap X})\xi,\xi> dg \\
                              &= \int f(g)<E_{g}\xi,\xi> dg. 
\end{align*}
Thus $\pi(F)=\int f(g)E_{g}$. 

Let $a \in Int(P)$ be given. Choose a sequence $(f_{n})$ as in Lemma \ref{crying out lemma} and Let $F_{n}(x):=\int f_{n}(g)1_{Xg}(x)$ for $x \in X$. Note that $F_{n}$ is uniformly bounded. By Lemma \ref{crying out lemma}, it follows that $F_{n}$ converges pointwise to $1_{X_{0}a}$. On the other hand, we have $\pi(F_{n})=\int f_{n}(g)E_{g}dg$. Since $g \to E_{g}$ is strongly continuous, it is easily verifiable that $\int f_{n}(g)E_{g}dg$ converges strongly to $E_{a}$. Hence $\pi(1_{X_{0}a})=E_{a}$. Clearly by definition $V_{a}^{*}\pi(\phi)V_{a}=\pi(R_{a}(\phi))$. This proves $(1)$ implies $(2)$.

Now assume $(2)$. The equality $V_{a}^{*}\pi(\phi)V_{a}=\pi(R_{a}(\phi))$ for $a \in P$ and $\phi \in C(X)$ translates to the fact that for $a \in P$ and $\xi \in \clh$, the push-forward measure $(\sigma_{a})_{*}(\mu_{\xi,\xi})=\mu_{V_{a}\xi,V_{a}\xi}$. Hence $V_{a}^{*}\pi(\phi)V_{a}=\pi(R_{a}(\phi))$ for $a \in P$ and $\phi \in B(X)$. Now by Remark \ref{covariance}, it is enough to show that $V_{a}\pi(\phi) V_{a}^{*}=\pi(R_{a^{-1}}(\phi))$ for $a \in Int(P)$ and $\phi \in B(X)$. Now let $a \in Int(P)$ and $\phi \in B(X)$ be given. Then by assumption $(2)$, we have $V_{a}^{*}\pi(R_{a^{-1}}\phi)V_{a}=\pi(\phi)$. Hence $V_{a}\pi(\phi)V_{a}^{*}=E_{a}\pi(R_{a^{-1}}(\phi))E_{a}$. By the strong continuity of $g \to E_{g}$, by assumption $(2)$ and Lemma \ref{decreasing sequence},  it follows that $\pi(1_{Xa})=E_{a}$.  Hence $V_{a}\pi(\phi)V_{a}^{*}=\pi(1_{Xa}R_{a^{-1}}(\phi))=\pi(R_{a^{-1}}(\phi))$. This completes the proof. \hfill $\Box$

\begin{thm}
\label{repn}
Let $X$ be a compact Hausdorff  space on which $P$ acts injectively. Let $\mathcal{G}:=X \rtimes P$. Assume that $\mathcal{G}$ has a Haar system. 
For $\phi \in C(X)$ and $f \in C_{c}(G)$, let $\phi \otimes f \in C_{c}(\mathcal{G})$ be defined by the equation $(\phi \otimes f)(x,g)=\phi(x)f(g)$. We denote $1\otimes f$ by $\widetilde{f}$.

Let $(\pi,V)$ be a covariant representation of $(X,P)$ on a Hilbert space $\clh$. Then there exists a representation $\lambda:C^{*}(\mathcal{G}) \to B(\clh)$ such that 
\begin{enumerate}
\item[(1)] For $f \in C_{c}(G)$,  $ \lambda(\widetilde{f})= \int \Delta(g)^{-\frac{1}{2}}f(g)W_{g^{-1}}dg$.  Here $\Delta$ is the modular function of the group $G$.
\item[(2)] For $\phi \in C(X)$ and $f \in C_{c}(G)$, $\lambda(\phi \otimes f)=\pi(\phi)\lambda(\widetilde{f})$. 
\end{enumerate}
\end{thm}

\textit{Proof of Theorem \ref{repn}.} Let $\phi \in C_{c}(\mathcal{G})$. We claim that  $G \ni g \to \pi(\phi_{g})W_{g^{-1}} \in B(\clh)$ is strongly continuous. Let $\Phi \in C_{c}(Y \rtimes G)$ be an extension of $\phi$.  Since $(\pi,V)$ is covariant, it follows that $E_{g}=\pi(1_{Xg \cap X})$. Now observe that $\pi(\widehat{\Phi_{g}})W_{g^{-1}}=\pi(\phi_{g})W_{g^{-1}}$. For $\pi(\widehat{\Phi_{g}})W_{g^{-1}}= \pi(\widehat{\Phi_{g}})E_{g^{-1}}W_{g^{-1}}= \pi(\widehat{\Phi_{g}}1_{Xg^{-1}})W_{g^{-1}}=\pi(\phi_{g})W_{g^{-1}}$. But $g \to \widehat{\Phi_{g}} \in C(X)$ is continuous. Hence $G \ni g \to \pi(\widehat{\Phi_{g}}) \in B(\clh)$ is strongly continuous and consequently $G \ni g \to \pi(\phi_{g})W_{g^{-1}} \in B(\clh)$ is strongly continuous.

For $\phi \in C_{c}(\mathcal{G})$, let $\lambda(\phi) \in B(\clh)$ be \[
                \lambda(\phi):= \int \Delta(g)^{-\frac{1}{2}}\pi(\phi_{g})W_{g^{-1}}dg.
              \]
 Also we have shown that if $\Phi \in C_{c}(Y \rtimes G)$ is an extension of $\phi \in C_{c}(\mathcal{G})$, then \[\lambda(\phi)=\int \Delta(g)^{-\frac{1}{2}}\pi_{Y}(\Phi_{g})W_{g^{-1}}dg.\]       
 For $\phi \in C_{c}(\mathcal{G})$, calculate as follows to find that
 \begin{align*}
 \lambda(\phi)^{*}&= \int W_{g} \pi(\overline{\phi_{g}}) \Delta(g)^{-\frac{1}{2}}dg \\
                           & = \int W_{g} \pi (\overline{\phi_{g}})W_{g}^{*}W_{g}\Delta(g)^{-\frac{1}{2}}dg \\
                            & = \int \pi(R_{g^{-1}}(\overline{\phi_{g}}))W_{g} \Delta(g)^{-\frac{1}{2}}dg \\
                             & = \int \pi(R_{g}(\overline{\phi_{g^{-1}}}))W_{g^{-1}}\Delta(g)^{\frac{1}{2}}\Delta(g)^{-1}dg \\
                             & = \int \pi((\phi^{*})_{g})W_{g^{-1}}\Delta(g)^{-\frac{1}{2}}dg \\
                             & = \lambda(\phi^{*}).
 \end{align*}             
 Thus $\lambda$ preserves the adjoint.
 
 Now let $\phi, \psi \in C_{c}(\mathcal{G})$ be given and let $\Phi, \Psi \in C_{c}(Y \rtimes G)$ be extensions of $\phi$ and $\psi$ respectively. Consider the function on $Y \rtimes G$ defined by the equation
 \[
 \Phi \circ \Psi(y,g)= \int \Phi (y,h)\Psi (y.h,h^{-1}g)1_{X}(y.h)dh.
 \]
 A simple application of the dominated convergence theorem together with the fact that $1_{X}(y.h)=1_{X_{0}}(y.h)$ a.e. for every $y$ implies that $\Phi \circ \Psi$ is continuous. Clearly $\Phi \circ \Psi $ is compactly supported and is an extension of $\phi * \psi$.
 
 Let $g \in G$ and $\xi \in \clh$ be given. Then 
 \begin{align*}
 <\pi_{Y}((\Phi \circ \Psi)_{g})\xi,\xi> &= \int \Phi \circ \Psi(x,g) d\mu_{\xi,\xi}(x)\\
                                                         & = \int \Big(\int  \Phi(x,h)\Psi(x.h,h^{-1}g)1_{X}(xh)dh\Big) d\mu_{\xi,\xi}(x) \\
                                                         & = \int \Big(\int \Phi(x,h)\Psi(x.h,h^{-1}g)1_{Xh^{-1}}(x)d\mu_{\xi,\xi}(x)\Big) dh \\
                                                         & = \int < \pi_{Y}(\Phi_{h})\pi_{Y}(R_{h}(\Psi_{h^{-1}g}))\pi_{Y}(1_{Xh^{-1}}) \xi, \xi> dh \\
                                                          & =  \int < \pi_{Y}(\Phi_{h})\pi_{Y}(R_{h}(\Psi_{h^{-1}g}))E_{h^{-1}} \xi, \xi> dh \\
  \end{align*}
 
 Thus  for $g \in G$, $ \pi_{Y}((\Phi \circ \Psi)_{g})= \int \pi_{Y}(\Phi_{h})\pi_{Y}(R_{h}(\Psi_{h^{-1}g}))E_{h^{-1}}dh$. Now calculate to find that
 \begin{align*}
 \lambda(\phi)\lambda(\psi) & = \int \Delta(gh)^{-\frac{1}{2}} \pi_{Y}(\Phi_{g})W_{g^{-1}}\pi_{Y}(\Psi_{h})W_{h^{-1}} dg dh \\
                                            & = \int  \Delta(gh)^{-\frac{1}{2}} \pi_{Y}(\Phi_{g})W_{g}^{*}\pi_{Y}(\Psi_{h})W_{g}W_{g^{-1}}W_{h^{-1}} dg dh \\
                                            & = \int   \Delta(gh)^{-\frac{1}{2}}  \pi_{Y}(\Phi_{g})\pi_{Y}(R_{g}(\Psi_{h}))W_{g^{-1}}W_{h^{-1}} dg dh \\
                                             & = \int   \Delta(gh)^{-\frac{1}{2}} \pi_{Y}(\Phi_{g})\pi_{Y}(R_{g}(\Psi_{h}))E_{g^{-1}}W_{h^{-1}g^{-1}}dg dh ~~\textrm{ [by Proposition \ref{commuting projections}]}\\
                                             & = \int \Big (  \int  \Delta(k)^{-\frac{1}{2}}\pi_{Y}(\Phi_{g})\pi_{Y}(R_{g}(\Psi_{g^{-1}k}))E_{g^{-1}}W_{k^{-1}}dk \Big) dg \\
                                              & = \int \Big ( \int \pi_{Y}(\Phi_{g})\pi_{Y}(R_{g}(\Psi_{g^{-1}k})E_{g^{-1}} dg \Big ) \Delta(k)^{-\frac{1}{2}}W_{k^{-1}}dk    \\
                                              & = \int \Delta(k)^{-\frac{1}{2}} \pi_{Y}((\Phi \circ \Psi)_{k})W_{k^{-1}} \\
                                              & = \lambda(\phi * \psi) . \\
                    \end{align*}
 Hence $\lambda$ preserves the multiplication. 

For $\phi \in B(X)$, one has $||\pi(\phi)|| \leq ||\phi||_{\infty}$ where $||.||_{\infty}$ is the sup norm on $B(X)$. Let $K$ be a compact subset of $G$. Then for $\phi \in C_{c}(\clg)$ with $supp(\phi) \subset X \times K$, observe that
\begin{align*}
 ||\lambda(\phi)|| &\leq \int \Delta(g)^{-\frac{1}{2}}||\pi(\phi_{g})||dg \\
                & \leq \int_{g \in K} \Delta(g)^{-\frac{1}{2}}||\phi_{g}|| dg \\
                & \leq \Big(\sup_{g \in K}\Delta(g)^{-\frac{1}{2}}\Big)||\phi||_{\infty} \int 1_{K}(g)dg.
\end{align*}

Thus it is clear that the map $\lambda:C_{c}(\mathcal{G}) \to B(\clh)$ is continuous when $C_{c}(\clg)$ is given the inductive limit topology and $B(\clh)$ is given the norm toplogy. By Renault's disintegration theorem, one obtains a bonafide representation $\lambda:C^{*}(\mathcal{G}) \to B(\clh)$. Conditions $(1)$ and $(2)$ follows just from definitions. This completes the proof. \hfill $\Box$

\section{The main theorem}

Let $V:P \to B(\clh)$ be an isometric representation with commuting range projections. Let $\mathcal{A}$ and $\Omega$ be as in sections 3-5. Denote the open set $\Omega Int(P)$ by $\Omega_{0}$.  Let $\pi: C(\Omega) \to B(\clh)$ be the representation induced by the inclusion $\mathcal{A} \subset B(\clh)$. Denote the extension to $B(\Omega)$ by $\pi$ itself. First let us show that $(\pi,V)$ is a covariant representation.



\begin{lmma}
\label{Covariant representation}
 The pair $(\pi,V)$ is a covariant representation.
\end{lmma}

\textit{Proof.} By definition, it follows that for $a \in P$, $V_{a}^{*}\pi(\phi)V_{a}=\pi(R_{a}(\phi))$ for $\phi \in C(X)$. Now fix $a \in P$. Choose a sequence $(f_{n}) \in C_{c}(G)$ as in Lemma \ref{crying out lemma}. Set $F_{n}:= \int f_{n}(g)E_{g}dg \in C(\Omega)$. Then by the strong continuity of $g \to E_{g}$, it is clear that $F_{n}$ converges strongly to $E_{a}$.
Now by definition, for $A \in \Omega$, $F_{n}(A)= \int f_{n}(g)1_{A}(g) dg$. By Proposition \ref{Haar system}, it follows that $F_{n}(A)= \int f_{n}(g)1_{\Omega}(Ag^{-1})dg$ for $A \in \Omega$. By Lemma \ref{crying out lemma}, it follows that $F_{n}$ converges pointwise to $1_{\Omega_{0}a}$. Hence $\pi(1_{\Omega_{0}a})=E_{a}$. The proof follows now from Proposition \ref{covariant}. \hfill $\Box$

\begin{ppsn}
\label{maintheorem1}
 Let $\clh$ be a Hilbert space and $V:P \to B(\clh)$ be an isometric representation with commuting range projections. Let $\Omega$ be as in Sections 3-5. Then there exists a $*$-homomorphism $\lambda:C^{*}(\Omega \rtimes P) \to B(\clh)$ such that for $f \in C_{c}(G)$, \[
                                                                                                                                                                                                                        \lambda(\widetilde{f})=\int \Delta(g)^{-\frac{1}{2}}f(g)W_{g^{-1}}dg.
                                                                                                                                                                                                                                         \]
Moreover the range of $\lambda$ is generated by $\{\int f(g)W_{g}dg: f \in C_{c}(G)\}$.
\end{ppsn}
\textit{Proof.}
For $f \in C_{c}(G)$, let $\widehat{f} \in C(\Omega)$ be defined by \[\widehat{f}(A):=\int f(g)1_{\Omega}(Ag)dg= \int f(g)1_{A}(g^{-1}).\]
By (4) of Remark \ref{Openness} and the fact that $\Omega \subset \Omega_{u}$, it follows that  $\{\widehat{f}:f \in C_{c}(G)\}$ separates points of $\Omega$. Thus by Proposition \ref{density}, the $*$-algebra generated by $\{\widetilde{f}:f \in C_{c}(G)\}$ is dense in $C^{*}(\Omega \rtimes P)$. Now the proof  follows directly from Lemma \ref{Covariant representation} and Proposition \ref{repn}. \hfill $\Box$

\begin{thm}
\label{maintheorem2}
 Let $\clh$ be a Hilbert space and $V:P \to B(\clh)$ be an isometric representation with commuting range projections. Let $\Omega_{u}:=\{A \in \mathcal{C}(G): P^{-1} \subset A ~\textrm{and~}P^{-1}A \subset A \}$ with the Vietoris topology. Consider the right action of  $P$  on $\Omega_{u}$ by right multiplication. Then there exists  $*$-homomorphism $\lambda:C^{*}(\Omega_{u} \rtimes P) \to B(\clh)$ such that for $f \in C_{c}(G)$, \[
                                                                                                                                                                                                                        \lambda(\widetilde{f})=\int \Delta(g)^{-\frac{1}{2}}f(g)W_{g^{-1}}dg.
                                                                                                                                                                                                                                         \]
Moreover the range of $\lambda$ is generated by $\{\int f(g)W_{g}dg: f \in C_{c}(G)\}$.
\end{thm}
\textit{Proof.} By Remark \ref{restriction}, it follows that $\Omega \rtimes P$ is isomorphic to the restriction $\Omega_{u} \rtimes P|_{\Omega}$ and $\Omega$ is an invariant subset of $\Omega_{u}$. Consider the natural map $res: C_{c}(\Omega_{u} \rtimes P) \to C_{c}(\Omega \rtimes P)$ which on $C_{c}(\Omega_{u} \rtimes P)$ is simply the restriction. Let $\widetilde{\lambda}: C^{*}(\Omega \rtimes P) \to B(\clh)$ be the representation as in Propositon \ref{maintheorem1}. Now one completes the proof by setting $\lambda:=\widetilde{\lambda} \circ res$. \hfill $\Box$

\begin{rmrk}
Proposition \ref{maintheorem2} says that the $C^{*}$-algebra of the groupoid $\Omega_{u} \rtimes P$ can be interpreted as  the 'universal' $C^{*}$-algebra which encodes the isometric representations with commuting range projections. However, the space $\Omega_{u}$ is quite large to describe explicitly even for the simple example of the quarter plane $[0,\infty)\times [0,\infty) \subset \mathbb{R}^{2}$. 
\end{rmrk}

We end this article by considering two well-known results which are a part of folklore in operator algebras.

\begin{xmpl}
 Let $P:=\mathbb{N}$ and $G:=\mathbb{Z}$ with the discrete topology. Consider the one-point compactification $\mathbb{N}_{\infty}:=\mathbb{N} \cup \{\infty\}$. The semigroup $\mathbb{N}$ acts on $\mathbb{N}_{\infty}$ by translation with the convention that $\infty + n = \infty$ for $n \in \mathbb{N}$. It is easy to verify that the map $\mathbb{N}_{\infty} \ni n \to (-\infty,n] \in \Omega_{u}$  is an $\mathbb{N}$-equivariant homeomorphism. Here $(-\infty,\infty]$ is just $\mathbb{N}$. The groupoid $\mathbb{N}_{\infty} \rtimes \mathbb{N}$ is amenable and  $C_{red}^{*}(\mathbb{N}_{\infty} \rtimes \mathbb{N})$ is just the Toeplitz-algebra. Now Theorem \ref{maintheorem2} is just the well-known Coburn's theorem.
\end{xmpl}

\begin{xmpl}Let $\mathbb{R}_{+}=[0,\infty)$.
 Let $P:=\mathbb{R}_{+}$ and $G:=\mathbb{R}$ with the usual Euclidean topology and addition as the group operation. Consider the one-point compactification $[0,\infty]:= [0,\infty) \cup \{\infty\}$. The semigroup $[0,\infty)$ acts on $\mathbb{R} \cup \{\infty\}$ by translation with the convention that $\infty + x=\infty$ for $x \in [0,\infty)$. It is easily verifiable that the map $[0,\infty] \ni x \to (-\infty,x] \in \Omega_{u}$ is a $\mathbb{R}_{+}$-equivariant homeomorphism. The groupoid $[0,\infty] \rtimes \mathbb{R}_{+}$ is amenable and $C_{red}^{*}([0,\infty] \rtimes [0,\infty))$ is the usual Wiener-Hopf algebra.
[See \cite{Renault_Muhly}] 

Observe that if $V:\mathbb{R}_{+} \to B(\clh)$ is an isometric representation then the range projections $\{E_{t}:=V_{t}V_{t}^{*}: t \geq 0 \}$ commutes. For if $t=r+s$ then $E_{t}E_{r}=V_{r}V_{s}V_{s}^{*}V_{r}^{*}V_{r}V_{r}^{*}=E_{t}$. Hence if $t >r$, then $E_{t}E_{r}=E_{t}$. Now the claim follows from the fact that $\mathbb{R}_{+}$ is totally ordered.  Thus if $V:\mathbb{R}_{+} \to B(\clh)$ is an isometric representation, then there exists a representation $\pi: \mathcal{W}([0,\infty),\mathbb{R}) \to B(\clh)$ such that 
\[
\pi(\widetilde{f})=\int_{0}^{\infty} f(t)V_{t} + \int_{-\infty}^{0}f(t)V_{t}^{*}
\]
for $f \in C_{c}(\mathbb{R})$.
\end{xmpl}

\bibliography{references}

\def\cprime{$'$} \def\cprime{$'$} \def\cprime{$'$}
\providecommand{\bysame}{\leavevmode\hbox to3em{\hrulefill}\thinspace}
\providecommand{\MR}{\relax\ifhmode\unskip\space\fi MR }
\providecommand{\MRhref}[2]{%
  \href{http://www.ams.org/mathscinet-getitem?mr=#1}{#2}
}
\providecommand{\href}[2]{#2}
\begin{thebibliography}{Mur96b}

\bibitem[Arv02]{Arv}
William Arveson, \emph{A short course on spectral theory}, Graduate Texts in
  Mathematics, vol. 209, Springer-Verlag, New York, 2002.

\bibitem[Cun08]{Cuntz}
Joachim Cuntz, \emph{{$C^*$}-algebras associated with the {$ax+b$}-semigroup
  over {$\Bbb N$}}, {$K$}-theory and noncommutative geometry, EMS Ser. Congr.
  Rep., 2008, pp.~201--215.

\bibitem[Exe03]{Exel_endo}
Ruy Exel, \emph{A new look at the crossed-product of a {$C^*$}-algebra by an
  endomorphism}, Ergodic Theory Dynam. Systems \textbf{23} (2003), no.~6,
  1733--1750.

\bibitem[HN95]{Hilgert_Neeb}
Joachim Hilgert and Karl-Hermann Neeb, \emph{Wiener-{H}opf operators on ordered
  homogeneous spaces. {I}}, J. Funct. Anal. \textbf{132} (1995), no.~1,
  86--118.

\bibitem[Li12]{Li-semigroup}
Xin Li, \emph{Semigroup {${\rm C}^*$}-algebras and amenability of semigroups},
  J. Funct. Anal. \textbf{262} (2012), no.~10, 4302--4340.

\bibitem[Li13]{Li13}
\bysame, \emph{Nuclearity of semigroup {$C^*$}-algebras and the connection to
  amenability}, Adv. Math. \textbf{244} (2013), 626--662.

\bibitem[MR82]{Renault_Muhly}
Paul~S. Muhly and Jean~N. Renault, \emph{{$C^{\ast} $}-algebras of
  multivariable {W}iener-{H}opf operators}, Trans. Amer. Math. Soc.
  \textbf{274} (1982), no.~1, 1--44.

\bibitem[Mur91]{Murphy91}
G.~J. Murphy, \emph{Ordered groups and crossed products of {$C^*$}-algebras},
  Pacific J. Math. \textbf{148} (1991), no.~2, 319--349.

\bibitem[Mur94]{Murphy94}
Gerard~J. Murphy, \emph{Crossed products of {$C^*$}-algebras by semigroups of
  automorphisms}, Proc. London Math. Soc. (3) \textbf{68} (1994), no.~2,
  423--448.

\bibitem[Mur96a]{Murphy}
\bysame, \emph{{$C^\ast$}-algebras generated by commuting isometries}, Rocky
  Mountain J. Math. \textbf{26} (1996), no.~1, 237--267.

\bibitem[Mur96b]{Murphy96}
\bysame, \emph{Crossed products of {$C^*$}-algebras by endomorphisms}, Integral
  Equations Operator Theory \textbf{24} (1996), no.~3, 298--319.

\bibitem[Nic87]{Nica_WienerHopf}
Alexandru Nica, \emph{Some remarks on the groupoid approach to {W}iener-{H}opf
  operators}, J. Operator Theory \textbf{18} (1987), no.~1, 163--198.

\bibitem[Nic90]{Nica90}
\bysame, \emph{Wiener-{H}opf operators on the positive semigroup of a
  {H}eisenberg group}, Linear operators in function spaces ({T}imi\c soara,
  1988), Oper. Theory Adv. Appl., vol.~43, Birkh\"auser, Basel, 1990,
  pp.~263--278.

\bibitem[RS15]{Jean_Sundar}
J~Renault and S~Sundar, \emph{Groupoids associated to {O}re semigroup actions},
  to appear in J. Operator Theory (2015), arXiv:1402.2762v2.

\end{thebibliography}
\bibliographystyle{amsalpha}

\noindent
{\sc S. Sundar}
(\texttt{sundarsobers@gmail.com})\\
         {\footnotesize  Chennai Mathematical Institute, H1 Sipcot IT Park, \\
Siruseri, Padur, 603103, Tamilnadu, INDIA.}

\end{document}